\documentclass[a4paper]{article}

\usepackage{amsmath}
\usepackage{amssymb}
\usepackage{mathrsfs}
\usepackage{amsthm}
\usepackage{tabularx}
\usepackage{enumerate} 
\usepackage{tikz-cd} 
\usepackage{hyperref} 
\usepackage{lscape} 

\newtheorem{thm}{Theorem}[section]
\newtheorem{defn}[thm]{Definition}
\newtheorem{lem}[thm]{Lemma}
\newtheorem{prop}[thm]{Proposition}
\newtheorem{coro}[thm]{Corollary}
\newtheorem*{rmk}{Remark}

\usepackage[top=1in, bottom=1.25in, left=1.25in, right=1.25in]{geometry}

\newcommand{\ml}[2]{\begin{tabular}{@{} >{$}#1<{$} @{}} #2 \end{tabular}}

\numberwithin{equation}{section}

\title{Graded $2$-generated axial algebras}
\author{Madeleine Whybrow\thanks{University of Primorska, madeleine.whybrow@famnit.upr.si}}
\date{}

\begin{document}

\maketitle

\begin{abstract}
    Axial algebras are non-associative algebras generated by semisimple idempotents whose adjoint actions obey a fusion law. Axial algebras that are generated by two such idempotents play a crucial role in the theory. We classify all primitive $2$-generated axial algebras whose fusion laws have two eigenvalues and all graded primitive $2$-generated axial algebras whose fusion laws have three eigenvalues. This represents a significant broadening in our understanding of axial algebras.
\end{abstract}

\section{Introduction}

Axial algebras are characterised as being non-associative algebras generated by semisimple idempotents, known as axes, whose adjoint actions obey a fusion law. The most important known examples of axial algebras are Jordan algebras and Majorana algebras, including the Griess algebra. The fusion laws that are associated with these two families of algebras, and their corresponding eigenvalues, have very specific forms. However, little is understood of the wider possibilities for fusion laws that admit non-trivial algebras. The results presented here offer a significant broadening of our understanding of axial algebras.

The \emph{Griess algebra} was constructed by Griess \cite{Griess82} as a 196,884 dimensional real commutative algebra whose automorphism group is equal to the Monster simple group. Conway \cite{Conway84} showed that there exists a bijection between the $2A$ conjugacy class of involutions in the Monster group and a generating set of idempotents in the Griess algebra, known as $2A$-\emph{axes}. Moreover, the adjoint action of these idempotents is semisimple and obeys the \emph{Monster fusion law}.

 Frenkel, Lepowsky and Meurman \cite{FLM88} introduced a class of infinite graded algebras known as \emph{vertex operator algebras}. Using this, they constructed the \emph{Moonshine module} that was later used by Borcherds \cite{Borcherds92} in his proof of Conway and Norton's Monstrous Moonshine conjectures. In general, the weight $2$ component of a vertex operator algebra has the structure of a real commutative algebra and is called a \emph{generalised Griess algebra}.

 Miyamoto \cite{Miyamoto96} showed that a generalised Griess algebra contains idempotents that are in bijection with certain involutions in its automorphism group, generalising Conway's observations concerning the $2A$-axes of the Griess algebra. Inspired by this work, and by a result of Sakuma \cite{Sakuma07}, Ivanov \cite{Ivanov09} introduced \emph{Majorana algebras} as an axiomatisation of certain properties of generalised Griess algebras.

Later, \emph{axial algebras} were introduced by Hall, Shpectorov and Rehren \cite{HRS15a} as a generalisation of Majorana algebras. Jordan algebras, and the closely related algebras of Jordan type, are also axial algebras and have been well studied in this context. There is also emerging evidence that axial algebras also play an important role in other areas of maths, such as in the study of cubic minimal cones, objects in algebraic geometry related to partial differential equations.

Jordan and Majorana algebras are of particular significance as their fusion laws admit a $C_2$-grading. In particular, this means that an algebra that obeys one of these fusion laws admits involutary algebra automorphisms that are given by the eigenspace decomposition of the adjoint actions of its axes. This gives these algebras a powerful and important link to groups generated by involutions. In the case of the Griess algebra, these involutions form the $2A$ conjugacy class of the Monster group. In this paper, we focus on fusion laws that are $G$-graded, for some group $G$.

The Griess algebra and Jordan algebras are known to be primitive - the one eigenspace of the adjoint action of a given axis is one dimensional, spanned by the axis itself. This seems to be a natural and useful additional assumption and is one which we take in this paper. To extend our methods to the non-primitive case would require more involved calculations and increased use of computational methods.

Given a generic fusion law, to classify the axial algebras that are generated by two axes is a crucial first step in understanding generic algebras that obey this law. Without a characterisation of the $2$-generated algebras of a fusion law, it is not possible to perform any constructive or classification results on algebras generated by larger numbers of axes.

Despite the importance of this question, relatively little is known of the $2$-generated algebras that obey a generic fusion law. Until now, the classification of $2$-generated axial algebras had been completed only for Majorana algebras \cite{IPSS10}, axial algebras of Monster type \cite{HRS15a} and axial algebras of Jordan type \cite{HRS15b}.

 We give a full classification of primitive axial algebras whose fusion laws have two eigenvalues and show that such an algebra is of dimension at most $2$. In the $3$-eigenvalue case, we restrict to classifying primitive algebras that are $G$-graded for some $G$. In this case, we give a full classification and show that the maximal dimension for such an algebra $3$.

These classification results also put strong restrictions on the possible fusion laws that can give rise to non-trivial algebras. In particular, we show that many fusion laws only give rise to non-trivial $2$-generated algebras for specific eigenvalues. In one such result, we consider a natural generalisation of the fusion law in the Jordan algebra case. We classify the possible eigenvalues for this law and show that the values that lead to Jordan algebras are in some sense distinguished.

The results presented here are constructive; in each case we give explicit values for the algebra product on a basis. This means that these classifications can be used in constructive results and programs, as has been successful in the case of Majorana algebras. In particular, the computational approaches of McInroy and Shpectorov \cite{MS18} and Pfeiffer and Whybrow \cite{PW18a} for constructing Majorana algebras and axial algebras of Monster type rely heavily on the classification of $2$-generated algebras in this case. These results will allow their methods to be expanded in order to construct $n$-generated algebras for the fusion laws considered here.

Certain low dimensional algebras will arise frequently in this work and we name them in line with, for example, \cite{IPSS10}. Note that if $V$ is an algebra and $X \subseteq V$ then the subalgebra of $V$ generated by $X$ is denoted $\langle \langle X \rangle \rangle$.

\begin{defn}
    Let $\mathbf{k}$ be a field and suppose that $V$ is a commutative $\mathbf{k}$-algebra. If $V = \langle \langle a_0 \rangle \rangle$ where $a_0$ is a non-zero idempotent (so that $V \cong \mathbf{k}$) then $V$ is referred to as the \emph{algebra} $1A$. If $V = \langle \langle a_0, a_1 \rangle \rangle$ where $a_0$ and $a_1$ are distinct non-zero idempotents such that $a_0a_1 = \alpha(a_0 + a_1)$ for some $\alpha \in \mathbf{k}$ then we say that $V$ is the \emph{algebra} $2B(\alpha)$.
\end{defn}

Note that the Norton-Sakuma dihedral algebra conventionally referred to as $2B$ is equal to $2B(0)$ in this notation.

The following are our main results.

\begin{thm}
    Let $V = \langle \langle a_0, a_1 \rangle \rangle$ be a primitive $(\mathcal{F}, *)$-axial algebra over $\mathbb{C}$ where $|\mathcal{F}| = 2$.
    Then $V$ is either the $1$-dimensional algebra $1A$ or the $2$-dimensional algebra $2B(\alpha)$ for some $\alpha \in \mathbb{C}$. If additionally $V$ is assumed to be graded then either $V$ is the algebra $1A$ or $\alpha \in \{-1, \frac{1}{2}\}$.
\end{thm}

\begin{thm}
    Let $V = \langle \langle a_0, a_1 \rangle \rangle$ be a graded primitive $(\mathcal{F}, *)$-axial algebra over $\mathbb{C}$ where $|\mathcal{F}| = 3$. Then either $V$ is the $1A$ algebra or $V$ is isomorphic to one of the algebras given in Table \ref{tab:mainthm} (on page \pageref{tab:mainthm}) where $\alpha, \beta, x \in \mathbb{C}$.
\end{thm}

The fusion laws that occur in the graded cases are all $C_2$-graded, with the exception of the fusion law given in Proposition \ref{prop:threeevalsfusionlaws} part (d), which admits a $C_3$-grading. There is opportunity for further work to determine the groups that can occur as Miyamoto groups of the algebras constructed in this paper.

In Section \ref{sec:axial_algebras}, we give all relevant definitions and preliminary results relating to fusion laws and axial algebras. We then classify the possible fusion laws that the algebras in question can obey and give these results in Section \ref{sec:fusion_laws}.

In Sections \ref{sec:2evals} and \ref{sec:3evals} we perform the constructions required for our main classification results, splitting this into the two and three eigenvalue cases. We then present the proofs of our main results in Section \ref{sec:main_results}.

\begin{rmk}
    In the three eigenvalue case, some calculations were performed with the aid of the computer algebra system Macaulay2 \cite{M2}. The code required to perform these calculations (described in Section \ref{sec:3evals}) can be found at \url{https://github.com/MWhybrow92/Graded-Algebras-Paper}.
\end{rmk}

\section{Preliminaries}
\label{sec:axial_algebras}

Henceforth, unless otherwise stated, all rings will be commutative and unital and all algebras will be commutative (but not necessarily associative).

\subsection{Fusion laws}
We use the definition of a fusion law that is introduced in a recent paper of De Medts, Peacock, Shpectorov and Van Couwenberghe \cite{DPSV19}. This definition is slightly more general than has been used previously.

\begin{defn}
    A \emph{fusion law} is a pair $(\mathcal{F}, *)$ where $\mathcal{F}$ is a set and $*:\mathcal{F} \times \mathcal{F} \rightarrow 2^{\mathcal{F}}$ is a map. A fusion law is called \emph{symmetric} if $x*y = y*x$ for all $x, y \in \mathcal{F}$.
\end{defn}

\begin{defn}
    Let $(\mathcal{F}, *)$ be a fusion law. An element $e \in \mathcal{F}$ is called a \emph{unit} if $e * x \subseteq \{x\}$ and $ x*e \subseteq \{x\}$ for all $x \in \mathcal{F}$.
\end{defn}

\begin{defn}
    Suppose that $(\mathcal{F},*_{\mathcal{F}})$ and $(\mathcal{H}, *_{\mathcal{H}})$ are two fusion laws. A \emph{morphism} from $(\mathcal{F},*_{\mathcal{F}})$ to $(\mathcal{H}, *_{\mathcal{H}})$ is a map $\xi : \mathcal{F} \rightarrow \mathcal{H}$ such that
    \[
        \xi(x *_{\mathcal{F}} y) \subseteq \xi (x) *_{\mathcal{H}} \xi (y)
    \]
    for all $x,y \in \mathcal{F}$. Here $\xi$ also denotes the extension of $\xi$ to a map $2^{\mathcal{F}} \rightarrow 2^{\mathcal{H}}$.
\end{defn}

This allows us to define a category $\mathbf{Fus}$ whose objects are fusion laws and whose morphisms are as given above.

\begin{defn}
    Let $(\mathcal{F},*_{\mathcal{F}})$ and $(\mathcal{H},*_\mathcal{H})$ be fusion laws. Then we say that $(\mathcal{H},*_\mathcal{H})$ is a \emph{sublaw} of $(\mathcal{F},*_{\mathcal{F}})$ if
    \begin{enumerate}[(i)]
        \item $\mathcal{H} \subseteq \mathcal{F}$;
        \item $x *_\mathcal{H} y \subseteq x *_\mathcal{F} y$ for all $x, y \in \mathcal{H}$.
    \end{enumerate}
\end{defn}

\subsection{Gradings}

The following is an important example of a fusion law.

\begin{defn}
    Let $G$ be a group. Then the pair $(G,*)$ where
    \begin{align*}
        *: G \times G &\rightarrow 2^G \\
        (g,h) &\mapsto \{gh\}
    \end{align*}
    is a \emph{group fusion law}.
\end{defn}

\begin{defn}
    Let $(\mathcal{F}, *)$ be a fusion law and let $(G,*)$ be a group fusion law. A $G$-grading of $(\mathcal{F}, *)$ is a morphism $f: (\mathcal{F}, *) \rightarrow (G,*)$. We say that a $G$-grading is \emph{adequate} if $f(\mathcal{F})$ generates $G$.
\end{defn}

Every fusion law admits a $G$-grading where $G$ is the trivial group; we call this the \emph{trivial} grading. All other gradings are considered \emph{non-trivial}. If a fusion law admits an adequate non-trivial grading then we say that it is \emph{graded}.

\begin{prop}[{\cite[Proposition 3.2]{DPSV19}}]
    \label{prop:mingrading}
    Every fusion law $(\mathcal{F}, *)$ admits a unique grading $f: \mathcal{F} \rightarrow \Gamma$ such that every grading of $(\mathcal{F}, *)$ factors through $(\Gamma, *)$. This grading is given by $f : \mathcal{F} \rightarrow \Gamma_{\mathcal{F}}$ where
    \[
        \Gamma_{\mathcal{F}} := \langle \gamma_x : x \in \mathcal{F} \mid \gamma_x \gamma_y = \gamma_z \textrm{ whenever } z \in x *y \rangle
    \]
    and
    \begin{align*}
        f: (\mathcal{F}, *) &\rightarrow (\Gamma_{\mathcal{F}}, *) \\
        x &\mapsto \gamma_x.
    \end{align*}
\end{prop}

\subsection{Axial algebras}

Let $(\mathcal{F}, *)$ be a symmetric fusion law and let $R$ be a unital commutative ring. Let $V$ be a commutative (not necessarily associative) $R$-algebra. Moreover, suppose that $\mathcal{F} \subseteq R$ and that $1_R \in \mathcal{F}$.

Throughout, if $a \in V$, $\lambda \in R$ and $\Lambda \subseteq R$, then we let $V^{a}_{\lambda} = \{ u \in V \mid a u  = \lambda u\}$ and $V^{a}_\Lambda = \sum_{\lambda \in \Lambda} V^{a}_{\lambda}$.

\begin{defn}
    \label{defn:axis}
    An idempotent $a \in V$ is an $(\mathcal{F}, *)$-\emph{axis} if
    \begin{enumerate}[(a)]
    \item $V = V_{\mathcal{F}}^a$;
    \item for all $\lambda, \mu \in \mathcal{F}$, $V_{\lambda}^{a}  V_{\mu}^{a} \subseteq V_{\lambda * \mu}^{a}$.
    \end{enumerate}
\end{defn}

\begin{defn}
    An $(\mathcal{F}, *)$-\emph{axial algebra} is a pair $(V,A)$ where $V$ is a commutative $R$-algebra and $A$ is a generating set of $(\mathcal{F}, *)$-axes. If $(\mathcal{F}, *)$ is graded (i.e. admits a non-trivial adequate grading) then we say that $(V,A)$ is \emph{graded}.
\end{defn}

We will later restrict to the case where $R$ is a commutative unital $\mathbf{k}$-algebra for some field $\mathbf{k}$ such that $\mathcal{F} \subseteq \mathbf{k}$. In this case, if $\lambda \neq \mu$ then $V^{a}_{\lambda} \cap V^{a}_{\mu} = 0$ and so $V = \bigoplus_{\lambda \in \mathcal{F}} V^{a}_\lambda$.

\begin{defn}
    An $(\mathcal{F}, *)$-axial algebra $(V,A)$ is $2$-\emph{generated} if there exists $a_0, a_1 \in A$ such that $V = \langle \langle a_0, a_1 \rangle \rangle$.
\end{defn}

Axial algebras generated by two axes are sometimes referred to as $\emph{dihedral}$. We reserve this term for the special case of $2$-generated axial algebras where the fusion law is $C_2$-graded (so that dihedral algebras admit dihedral Miyamoto groups).

\begin{defn}
    An $(\mathcal{F}, *)$-axis $a \in V$ is \emph{primitive} if $V_1^{a} = Ra$. An $(\mathcal{F}, *)$-axial algebra $(V,A)$ is \emph{primitive} if $a$ is primitive for all $a \in A$.
\end{defn}

The following definition is introduced in \cite{FMS20}.

\begin{defn}
    A vector $v \in V$ is \emph{free} if its annihilator ideal $\mathrm{Ann}_R(v) := \{r \in R \mid rv = 0 \}$ is trivial.
\end{defn}

If $R$ is a field then all non-zero vectors are free.

\subsection{Miyamoto groups}

In the following, let $R$ be a commutative unital $\mathbf{k}$-algebra (so that $R$ can be considered as a ring containing a copy of $\mathbf{k}$ such that $1_R= 1_{\mathbf{k}}$). Let $(\mathcal{F}, *)$ be a symmetric fusion law such that $\mathcal{F} \subseteq \mathbf{k}$ and $1_R \in \mathcal{F}$. We will show how gradings of fusion laws lead to automorphisms of axial algebras. This provides the crucial link between axial algebras and group theory and motivates our focus on the graded case.

\begin{defn}
    Let $R^\times$ be the group of invertible elements of the $\mathbf{k}$-algebra $R$. An $R$-character of $G$ is a group homomorphism $\chi: G \rightarrow R^\times$. We let $\chi_R(G)$ denote the $R$-\emph{character group} of $G$, i.e. all $R$-characters of $G$, with group operation given by multiplication in $R^\times$.
\end{defn}

\begin{defn}
    Let $V$ be a commutative $R$-algebra and let $a \in V$ be an $(\mathcal{F}, *)$-axis. Suppose that $(\mathcal{F}, *)$ admits a $G$-grading $f$, for some group $G$. Then, for each $\chi \in \chi_R(G)$, we define a linear map such that
    \begin{align*}
        \tau_{a, \chi}: V &\rightarrow V \\
                        v &\mapsto \chi(f(\lambda)) v
    \end{align*}
    for all $v \in V_\lambda^{a}$. If $(V,A)$ is an axial algebra then we define its \emph{Miyamoto group} to be
    \[
        \mathrm{Miy}(V,A) := \langle \tau_{a, \chi} \mid a \in A, \chi \in \chi_R(G) \rangle.
    \]
\end{defn}

Note that, as $\mathcal{F} \subseteq \mathbf{k}$, the spaces $V_\lambda^{a}$ for $\lambda \in \mathcal{F}$ intersect trivially and so the maps $\tau_{a, \chi}$ are well-defined.

We are often interested in the case where $(\mathcal{F}, *)$ admits a $C_2$-grading and $R$ a field of characteristic $0$. In this case, $\chi_R(G) = \{1, \chi\}$, where $\chi$ maps the non-trivial element of $G$ to $-1 \in R$, and $\tau_{a, \chi}$ is the Miyamoto involution, as defined in \cite[Proposition 3.4]{HRS15a}.

\begin{prop}
    Let $V$ be a commutative $R$-algebra and let $a \in V$ be an $(\mathcal{F}, *)$-axis. If $\chi \in \chi_R(G)$ then $\tau_{a, \chi}$ is an automorphism of $V$.
\end{prop}

\subsection{Universal primitive axial algebras}
\label{sec:universalprimitive}

Now let $\mathbf{k}$ be a field and let $(\mathcal{F}, *)$ be a symmetric fusion law such that $\mathcal{F} \subseteq \mathbf{k}$ and $1_{\mathbf{k}} \in \mathcal{F}$.

In this section we construct the universal primitive $n$-generated $(\mathcal{F}, *)$-axial algebra; an idea that will be key in the proof of our main results. Given a fusion law, we want to set up an algebra $V$ so that every primitive $(\mathcal{F}, *)$-axial algebra generated by $n$ (or fewer) free axes occurs as a quotient of $V$. These sections closely follow the method given in \cite[Section 4]{HRS15a}. We outline the construction of the universal primitive algebra but do not provide any proofs; these can be easily adapted from their counterparts in \cite{HRS15a}.

The authors of \cite{HRS15a} consider axial algebras that admit an associating bilinear form $\langle \, , \, \rangle: V \times V \rightarrow \mathbf{k}$; we consider the slightly more general case of axial algebras that admit a linear map $\varphi_a: V \rightarrow \mathbf{k}$ for each axis $a$ of $V$. This approach was also taken in \cite{Rehren15}.

This section is based on the following observation.

\begin{prop}
    Let $R$ be a commutative unital $\mathbf{k}$-algebra and let $V$ be a commutative $R$-algebra. Suppose that $a \in V$ is a free $(\mathcal{F}, *)$-axis. Then $a$ is primitive if and only if there exists a linear map $\varphi_a : V \rightarrow R$ such that $v = \varphi_a(v)a$ for all $v \in V_1^a$.
\end{prop}

\begin{proof}
    As $a$ is an $(\mathcal{F}, *)$-axis and $\mathcal{F} \subseteq \mathbf{k}$, $V = \bigoplus_{\lambda \in \mathcal{F}}V_{\lambda}^{a}$. If $a$ is primitive then, as it is free, the projection $\varphi_a: V \rightarrow V_1^a = Ra$ is a well-defined linear map. Moreover, $V_1^a = Ra$ so if $v \in V_1^a$ then $v = \varphi_a(v)a$, as required. The converse is clear.
\end{proof}

We first consider the category $\mathcal{C}^{P}_0$ constructed as follows. Its objects consist of tuples $(R, V, A)$ such that
\begin{enumerate}[(i)]
    \item $R$ is a commutative unital $\mathbf{k}$-algebra;
    \item $V$ is a commutative $R$-algebra;
    \item $A = (a_0, a_1, \ldots, a_{n-1})$ is an ordered tuple of free elements such that $V = \langle \langle A \rangle \rangle$ and such that $V$ is equipped with a linear map $\varphi_{a} : V \rightarrow R$ for each $a \in A$.
\end{enumerate}
We call $R$ the \emph{coefficient ring} of $V$.

A morphism from $(R, V, A)$ to $(R', V', A')$ is a pair $(\phi, \psi)$ such that
\begin{enumerate}[(i)]
    \item $\phi$ is an algebra homomorphism from $R$ to $R'$ such that $\phi$ restricted to the field $\mathbf{k}$ is the identity;
    \item $\psi$ is an algebra homomorphism from $V$ to $V'$;
    \item $\psi(rm) = \phi(r)\psi(m)$ for all $r \in R$ and $m \in V$;
    \item $\psi(a_i) = a'_i$ for $i < n$, where $A = (a_0, a_1, \ldots, a_{n-1})$ and $A' = (a'_0, a'_1, \ldots, a'_{n-1})$;
    \item $ \phi(\varphi_{a}(m)) = \varphi_{\psi(a)}(\psi(m))$ for all $a \in A$ and $m \in V$.
\end{enumerate}

We will now construct the initial object, $\hat{V}$, of $\mathcal{C}^{P}_0$. Let $\hat{M}$ be the free commutative, non-associative magma with marked generators $\hat{A} :=(\hat{a}_0, \hat{a}_1, \ldots, \hat{a}_{n-1})$. Let
\[
    \hat{R} := \mathbf{k}[\{\varphi_{\hat{a}}(\hat{x})\}_{\hat{a} \in \hat{A}, \hat{x} \in \hat{M}}]
\]
and let $\hat{V} = \hat{R}\hat{M}$, the set of all formal linear combinations $\sum_{\hat{x} \in \hat{M}} \alpha_{\hat{x}}\hat{x}$ where $\alpha_{\hat{x}} \in \hat{R}$. Then $\hat{V}$ is a commutative $\hat{R}$-algebra where the algebra product is defined using the operation in $\hat{M}$. For $\hat{a} \in \hat{A}$, we can define the linear map $\varphi_{\hat{a}}$ via
\[
    \varphi_{\hat{a}}: \sum_{\hat{x} \in \hat{M}} \alpha_{\hat{x}}\hat{x}\mapsto \sum_{\hat{x} \in \hat{M}} \alpha_{\hat{x}} \varphi_{\hat{a}}(\hat{x}).
\]
So $(\hat{R}, \hat{V}, \hat{A})$ is an object of $\mathcal{C}^{P}_0$.

\begin{prop}
    \label{prop:4.2}
    The algebra $(\hat{R}, \hat{V}, \hat{A})$ is the initial object of $\mathcal{C}^{P}_0$. That is to say, for every object $(R, V, A)$ of $\mathcal{C}^{P}_0$, there exists a unique morphism $(\psi_V, \phi_V)$ from $(\hat{R}, \hat{V}, \hat{A})$ to $(R, V, A)$.
\end{prop}

We say that an object $(R, V, A)$ of $\mathcal{C}^{P}_0$ is \emph{rigid} if both $\psi_V$ and $\phi_V$ are surjective. That is to say, $V$ is generated as an algebra by $A$ and $R$ is generated by the field  and the values of the maps $\varphi_{a}$ for $a \in A$ on the submagma $M$ of $V$ generated by $A$. We let $\mathcal{C}^{P}$ be the subcategory of $\mathcal{C}^{P}_0$ whose objects are the isomorphism class representatives of rigid objects of $\mathcal{C}^{P}_0$.

\begin{prop}
    Every algebra $V$ in $\mathcal{C}^{P}$ is the quotient, up to isomorphism, of $\hat{V}$ over the ideal $I_V = \ker \psi_V$ and its coefficient ring $R$ is the quotient, up to isomorphism, of $\hat{R}$ over the ideal $J_V = \ker \phi_V$.
\end{prop}

\begin{defn}
    If $I \subseteq \hat{V}$ and $J \subseteq \hat{R}$ are ideals then we say that $I$ and $J$ \emph{match} if $I = I_V$ and $J = J_V$ for some $V$ from $\mathcal{C}^{P}$.
\end{defn}

\begin{prop}
    Suppose that $I \subseteq \hat{V}$ and $J \subseteq \hat{R}$ are ideals. Then $I$ and $J$ match if and only if $J \neq \hat{R}$ and, additionally, $J\hat{V} \subseteq I$ and $\varphi_{\hat{a}}(I) \subseteq J$ for all $\hat{a} \in \hat{A}$.
\end{prop}

Choose $X \leq \hat{V}$. Then we let $\mathcal{C}^{P}_X$ be the subcategory of $\mathcal{C}^{P}$ consisting of all algebras $V$ from $\mathcal{C}^{P}$ such that $X \subseteq I_V$. Let
\begin{enumerate}[(i)]
    \item $I_0$ be the ideal of $\hat{V}$ generated by $X$;
    \item $J_X$ be the ideal of $\hat{R}$ generated by the sets
    \[
    \{\varphi_{\hat{a}}(i) \mid \hat{a} \in \hat{A}, i \in I_0 \} \textrm{ and } \{ r \mid \exists a \in \hat{A} \textrm{ s.t. } ra \in I_0 \};
    \]
    \item $I_X$ be the ideal of $\hat{V}$ generated by $X$ and $J_X\hat{V}$.
\end{enumerate}

\begin{prop}
    \label{prop:cpx}
    The subcategory $\mathcal{C}^{P}_X$ is non-empty if and only if $J_X \neq \hat{R}$. Furthermore, if the latter holds then $I_X$ and $J_X$ match and hence correspond to some algebra $V_X$ of $\mathcal{C}^{P}$. The algebra $V_X$ is the initial object of $\mathcal{C}^{P}_X$ and $\mathcal{C}^{P}_X$ consists of all quotients of $V_X$.
\end{prop}

We say that $X$ is a \emph{presentation} for $V_X$ and that the elements of $X$ are the \emph{algebra relators} of $V_X$.

\begin{defn}
    Suppose that $\Lambda \subseteq R$. Then we let $f_{\Lambda} = \prod_{\lambda \in \Lambda} (x - \lambda) \in R[x]$.
\end{defn}

\begin{thm}
    \label{thm:universalprimitive}
    An algebra $(R, V, A)$ from $\mathcal{C}^{P}$ is a primitive $(\mathcal{F}, *)$-axial algebra generated by free axes if and only if each of the following conditions holds
    \begin{enumerate}[(i)]
        \item $\hat{a}\hat{a} - \hat{a} \in I_V$ for all $\hat{a} \in \hat{A}$;
        \item $f_{\mathcal{F} \backslash \{1\}}(ad_{\hat{a}})(\hat{x} - \varphi_{\hat{a}}(\hat{x})\hat{a}) \in I_V$ for all $\hat{a} \in \hat{A}$ and $\hat{x} \in \hat{M}$;
        \item $f_{\lambda*\mu}(ad_{\hat{a}})( (f_{\mathcal{F}\backslash\{\lambda\}}(ad_{\hat{a}})(\hat{x})) (f_{\mathcal{F}\backslash\{\mu\}}(ad_{\hat{a}})(\hat{y})) ) \in I_V $ for all $\lambda, \mu \in \mathcal{F}$, $\hat{a} \in \hat{A}$ and $\hat{x}, \hat{y} \in \hat{M}$.
    \end{enumerate}
    If $X$ consists of the vectors listed above then the category of $n$-generated primitive $(\mathcal{F}, *)$-axial algebras coincides with $\mathcal{C}^{P}_{X}$ and $V_X$ is the universal $n$-generated primitive $(\mathcal{F}, *)$-axial algebra, of which all others are quotients.
\end{thm}

The following are useful results concerning primitive axial algebras.

\begin{prop}
    \label{prop:conditionii}
    Suppose that $(R, V, A)$ is an object of $\mathcal{C}^{P}$ that obeys conditions (i) and (ii) above. Then
    \[
        f_{\mathcal{F}}(ad_{\hat{a}})(\hat{x}) \in I_V
    \]
    for all  $\hat{a} \in \hat{A}$ and $\hat{x} \in \hat{M}$. Equivalently, $f_{\mathcal{F}}(ad_a)(x) = 0$ for all $a \in A$ and all $x \in V$.
\end{prop}

\begin{proof}
    First, condition (ii) above says that $f_{\mathcal{F} \backslash \{1\}}(ad_{\hat{a}})(\hat{x} - \varphi_{\hat{a}}(\hat{x})\hat{a}) \in I_V$ for all $\hat{a} \in \hat{A}$ and $\hat{x} \in \hat{M}$. Then
    \begin{align*}
            f_{\mathcal{F}}(ad_{\hat{a}})(\hat{x})+ I_V &= (ad_{\hat{a}} - 1)f_{\mathcal{F}\backslash\{1\}}(ad_{\hat{a}})(\hat{x}) + I_V \\
            &= (ad_{\hat{a}} - 1)f_{\mathcal{F}\backslash\{1\}}(ad_{\hat{a}})(\varphi_{\hat{a}}(\hat{x})\hat{a}) + I_V \\
            &= \varphi_{\hat{a}}(\hat{x}) f_{\mathcal{F}\backslash\{1\}}(ad_{\hat{a}})(ad_{\hat{a}} - 1)(\hat{a}) + I_V \\
            &= \varphi_{\hat{a}}(\hat{x}) f_{\mathcal{F}\backslash\{1\}}(ad_{\hat{a}})(\hat{a}\hat{a} - \hat{a}) + I_V \\
            &= I_V.
    \end{align*}
    So that $f_{\mathcal{F}}(ad_{\hat{a}})(\hat{x}) \in I_V$ for all $\hat{a} \in \hat{A}$ and $\hat{x} \in \hat{M}$.

    Let $M$ be the submagma in $V$ generated by $A$. Then $f_{\mathcal{F}}(ad_{\hat{a}})(\hat{x}) \in I_V$ for all $\hat{a} \in \hat{A}$ and $\hat{x} \in \hat{M}$ if and only if $f_{\mathcal{F}}(ad_a)(x) \in I_V$ for all $a \in A$ and $x \in M$. The final claim then follows from the fact that $M$ spans $V$.
\end{proof}

\begin{prop}
    \label{prop:projections}
    Suppose that $(V,A)$ is a primitive $(\mathcal{F}, *)$-axial algebra. Then for all $a \in A$
    \begin{enumerate}[(i)]
        \item $\varphi_a(a) = 1$;
        \item $\varphi_a(v_\lambda) = 0$ for all $v_\lambda \in V_\lambda^a$ and for all $\lambda \in \mathcal{F} \backslash \{1\}$.
    \end{enumerate}
\end{prop}

\begin{proof}
    As $aa - a = 0$, condition (ii) of Theorem \ref{thm:universalprimitive} implies that
    \[
        f_{\mathcal{F} \backslash \{1\}}(ad_{a})(a - \varphi_{a}(a)a) =  \prod_{\lambda \in \mathcal{F} \backslash \{1\}} (1  - \lambda) (a - \varphi_{a}(a)a) = 0
    \]
    As $\prod_{\lambda \in \mathcal{F} \backslash \{1\}} (1  - \lambda) \neq 0$ is invertible and as $a$ is a non-zero vector, we conclude that $\varphi_{a}(a) = 1$.

    Now, as $v_\lambda \in V_\lambda^a$, $(ad_a - \lambda)(v_\lambda) = 0$. Then
    \begin{align*}
        f_{\mathcal{F} \backslash \{1\}}(ad_{a})(v_\lambda - \varphi_{a}(v_\lambda)a) &=  f_{\mathcal{F} \backslash \{1, \lambda\}}(ad_{a})(ad_a - \lambda)(v_\lambda - \varphi_{a}(v_\lambda)a) \\
        & = f_{\mathcal{F} \backslash \{1, \lambda\}}(ad_{a}) (\lambda - 1)(\varphi_{a}(v_\lambda)a) \\
        & = - \varphi_{a}(v_\lambda) \prod_{\mu \in \mathcal{F} \backslash \{1\}} ( 1 - \mu) a \\
        &= 0.
    \end{align*}
    So, as before, we can conclude that $\varphi_{a}(v_\lambda) = 0$.
\end{proof}

\subsection{Semiautomorphisms of primitive axial algebras}

As before, let $\mathbf{k}$ be a field and let $(\mathcal{F}, *)$ be a symmetric fusion law such that $\mathcal{F} \subseteq \mathbf{k}$ and $1_{\mathbf{k}} \in \mathcal{F}$. Let $(R, V, A)$ be the universal $2$-generated primitive $(\mathcal{F}, *)$-axial algebra with $A = (a_0, a_1)$.

\begin{defn}
    Let $\phi$ be an automorphism of $R$ that fixes $\mathbf{k}$ element-wise. If $\psi$ is a bijection $V \rightarrow V$ preserving addition and multiplication such that
    \[
        (rm)^\psi = r^\phi m^\psi \textrm{ and } \varphi_{a^\psi}(m^\psi) = \varphi_{a}(m)^\phi
    \]
    then we say that $\psi$ is a $\phi$-\emph{automorphism} of $V$. In general, we call $\phi$-automorphisms of $V$ \emph{semiautomorphisms}.
\end{defn}

\begin{prop}
    \label{prop:sigma}
    The permutation $(a_0, a_1) \in \mathrm{Sym}(A)$ uniquely extends to a semiautomorphism $\sigma$ of $V$.
\end{prop}

\begin{proof}
    We first consider the algebra $\hat{V}$ and define the permutation $\hat{\sigma} \in \mathrm{Sym}(\hat{A})$ to be such that
    \[
        \hat{\sigma} = (\hat{a}_0, \hat{a}_1).
    \]
    Then $\hat{\sigma}$ clearly extends uniquely to an automorphism $\hat{\psi}_\sigma$ of the free commutative magma $\hat{M}$. Further, it extends to  permutation $\hat{\phi}_\sigma$ of $\hat{R}$ which acts on the indeterminates of $\hat{R}$ via
    \[
        \varphi_{a_0}(\hat{x})^{\hat{\phi}_\sigma} = \varphi_{a_1}(\hat{x})
    \]
    for all $\hat{x} \in \hat{M}$. This is clearly an automorphism of $\hat{R}$ that fixes the field $\mathbf{k}$. We can now see that $\hat{\psi}_\sigma$ uniquely extends to $\hat{V}$ via
    \[
    \left( \sum_{\hat{x} \in \hat{M}} \alpha_{\hat{x}} \hat{x} \right)^{\hat{\psi}_{\sigma}} = \sum_{\hat{x} \in \hat{M}} \alpha_{\hat{x}}^{\hat{\phi}_{\sigma}} \hat{x}^{\hat{\psi}_{\sigma}}.
    \]
    Crucially, the maps $\hat{\phi}_\sigma$ and $\hat{\psi}_\sigma$ preserve the set $X$ of algebra relators given in Theorem \ref{thm:universalprimitive} and therefore also the ideals $I_X$ and $J_X$. In particular, these maps descend to automorphisms $\phi_\sigma$ and $\psi_\sigma$ of $R$ and $V$ respectively and we have an automorphism $\sigma := (\phi_\sigma, \psi_\sigma) $.
\end{proof}

\begin{prop}
    \label{prop:isomorphicfusionlaws}
    Suppose that $(\mathcal{F}, *)$ and $(\mathcal{H}, *)$ are two fusion laws such that $\mathcal{F}, \mathcal{H} \subseteq \mathbf{k}$, $1_{\mathbf{k}} \in \mathcal{F}$ and $1_{\mathbf{k}} \in \mathcal{H}$. Suppose further that there exists an isomorphism $\rho: (\mathcal{F}, *) \rightarrow (\mathcal{H}, *)$ such that $\rho(\mathcal{F}) = \mathcal{H}$ and $\rho(1_{\mathbf{k}}) = 1_{\mathbf{k}}$. Then their corresponding universal primitive $n$-generated axial algebras are equal.
\end{prop}

\begin{proof}
    We see that the only point in the construction of the universal algebra that involves the fusion law is in the algebra relators defined in Theorem \ref{thm:universalprimitive}. These are all preserved by such an isomorphism of fusion laws and so the resulting algebras must be equal.
\end{proof}

\subsection{Sublaws and universal algebras}

Again, let $\mathbf{k}$ be a field and let $(\mathcal{F}, *)$ be a symmetric fusion law such that $\mathcal{F} \subseteq \mathbf{k}$ and $1_{\mathbf{k}} \in \mathcal{F}$.

\begin{thm}
    \label{thm:sublaw}
    Let $(\mathcal{H}, *_{\mathcal{H}})$ be a sublaw of $(\mathcal{F}, *_{\mathcal{F}})$. Then the universal $n$-generated primitive $(\mathcal{H}, *_{\mathcal{H}})$-axial algebra occurs as a quotient of the universal $n$-generated primitive $(\mathcal{F}, *_{\mathcal{F}})$-axial algebra.
\end{thm}

\begin{proof}
    Let the sets $X_{\mathcal{F}}, X_{\mathcal{H}} \subseteq \hat{V}$ consist of the algebra relators described in Theorem \ref{thm:universalprimitive} for the fusion laws $(\mathcal{F}, *_{\mathcal{F}})$ and $(\mathcal{H}, *_{\mathcal{H}})$ respectively. Let $I_{0, \mathcal{F}}$ and $I_{0, \mathcal{H}}$ be the ideals generated by $X_{\mathcal{F}}$ and $X_{\mathcal{H}}$ respectively. We will show that $I_{0, \mathcal{F}} \subseteq I_{0, \mathcal{H}}$.

    The requirement (i) of Theorem \ref{thm:universalprimitive} says that $\hat{a}\hat{a} - \hat{a} \in X_{\mathcal{F}}$ for all $\hat{a} \in \hat{A}$. Clearly also $\hat{a}\hat{a} - \hat{a} \in X_{\mathcal{H}}$.

    From requirement (ii), $f_{\mathcal{H} \backslash \{1\}}(ad_{\hat{a}})(\hat{x} - \varphi_{\hat{a}}(\hat{x})\hat{a}) \in X_{\mathcal{H}}$ for all $\hat{a} \in \hat{A}$ and $x \in \hat{M}$. Thus
    \[
        f_{\mathcal{F} \backslash \{1\}}(ad_{\hat{a}})(\hat{x} - \varphi_{\hat{a}}(\hat{x})\hat{a}) = f_{\mathcal{F} \backslash \mathcal{H}}(ad_{\hat{a}}) f_{\mathcal{H} \backslash \{1\}}(ad_{\hat{a}})(\hat{x} - \varphi_{\hat{a}}(\hat{x})\hat{a}) \in I_{0, \mathcal{H}}.
    \]

    We now turn our attention to the vectors given by requirement (iii). From Proposition \ref{prop:conditionii}, $f_{\mathcal{H}}(ad_{\hat{a}})(\hat{x}) \in I_{0, \mathcal{H}}$ for all $\hat{a} \in \hat{A}$ and $\hat{x} \in \hat{M}$. Let $\lambda, \mu \in \mathcal{F}$ and assume that at least one of $\lambda$ and $\mu$ does not lie in $\mathcal{H}$. Without loss of generality, we assume that $\lambda \notin \mathcal{H}$. Then $f_{\mathcal{F}\backslash \{\lambda\}}(ad_{\hat{a}})(\hat{x}) = f_{(\mathcal{F} \backslash \mathcal{H})\backslash \{\lambda\}}(ad_{\hat{a}})f_{\mathcal{H}}(ad_{\hat{a}})(\hat{x}) \in I_{0, \mathcal{H}}$ and so
    \begin{align*}
        f_{\lambda*_{\mathcal{F}}\mu}(ad_{\hat{a}})( (f_{\mathcal{F}\backslash\{\lambda\}}(ad_{\hat{a}})(\hat{x})) (f_{\mathcal{F}\backslash\{\mu\}}(ad_{\hat{a}})(\hat{y})) ) \in I_\mathcal{0, H}
    \end{align*}
    for all $\hat{x}, \hat{y} \in \hat{M}$ as required.

    We now need to consider the case where $\lambda, \mu \in \mathcal{H}$. If $\lambda \in \mathcal{H}$ then, as $f_{\mathcal{H}}(ad_{\hat{a}})(\hat{x}) = (ad_{\hat{a}} - \lambda)f_{\mathcal{H}\backslash \{\lambda\}}(ad_{\hat{a}})(\hat{x}) \in I_{0, \mathcal{H}}$,  we have
    \[
        ad_{\hat{a}} (f_{\mathcal{H}\backslash \{\lambda\}}(ad_{\hat{a}})(\hat{x})) + I_{0, \mathcal{H}} = \lambda f_{\mathcal{H}\backslash \{\lambda\}}(ad_{\hat{a}})(\hat{x}) + I_{0, \mathcal{H}}.
    \]
    Thus
    \begin{align*}
        f_{\mathcal{F}\backslash \{\lambda\}}(ad_{\hat{a}})(\hat{x}) + I_{0, \mathcal{H}} &= f_{\mathcal{F}\backslash \mathcal{H}}(ad_{\hat{a}})f_{\mathcal{H}\backslash \{\lambda\}}(ad_{\hat{a}})(\hat{x}) + I_{0, \mathcal{H}} \\
        &= \prod_{\nu \in \mathcal{F}\backslash \mathcal{H}}(ad_{\hat{a}} - \nu)f_{\mathcal{H}\backslash \{\lambda\}}(ad_{\hat{a}})(\hat{x}) + I_{0, \mathcal{H}} \\
        &= \prod_{\nu \in \mathcal{F}\backslash \mathcal{H}}(\lambda - \nu)f_{\mathcal{H}\backslash \{\lambda\}}(ad_{\hat{a}})(\hat{x}) + I_{0, \mathcal{H}}. \\
    \end{align*}
    Thus, modulo $I_{0, \mathcal{H}}$, $f_{\mathcal{F}\backslash \{\lambda\}}(ad_{\hat{a}})(\hat{x})$ is a scalar multiple of $f_{\mathcal{H}\backslash \{\lambda\}}(ad_{\hat{a}})(\hat{x})$

    Finally, as $(\mathcal{H}, *_{\mathcal{H}})$ is a sublaw of $(\mathcal{F}, *_{\mathcal{F}})$, $\lambda *_{\mathcal{H}} \mu \subseteq \lambda *_{\mathcal{F}} \mu$ and so for all $\hat{x}, \hat{y} \in \hat{M}$
    \begin{align*}
        f_{\lambda*_{\mathcal{H}} \mu}(ad_{\hat{a}})( (f_{\mathcal{H}\backslash\{\lambda\}}(ad_{\hat{a}})(\hat{x})) (f_{\mathcal{H}\backslash\{\mu\}}(ad_{\hat{a}})(\hat{y})) ) &\in I_{0, \mathcal{H}} \\
        \Rightarrow f_{\lambda*_{\mathcal{F}} \mu}(ad_{\hat{a}})( (f_{\mathcal{H}\backslash\{\lambda\}}(ad_{\hat{a}})(\hat{x})) (f_{\mathcal{H}\backslash\{\mu\}}(ad_{\hat{a}})(\hat{y})) ) &\in I_{0, \mathcal{H}} \\
        \Rightarrow f_{\lambda*_{\mathcal{F}} \mu}(ad_{\hat{a}})( (f_{\mathcal{F}\backslash\{\lambda\}}(ad_{\hat{a}})(\hat{x})) (f_{\mathcal{F}\backslash\{\mu\}}(ad_{\hat{a}})(\hat{y})) )& \in I_{0, \mathcal{H}}
    \end{align*}
    as required and so $I_{0, \mathcal{F}} \subset I_{0, \mathcal{H}}$.

    As $I_{0, \mathcal{F}} \subset I_{0, \mathcal{H}}$, $I_{X_\mathcal{F}} \subseteq I_{X_\mathcal{H}}$ and $J_{X_\mathcal{F}} \subseteq J_{X_\mathcal{H}}$ and so the universal $n$-generated primitive $(\mathcal{H}, *)$-axial algebra occurs as a quotient of the universal $n$-generated primitive $(\mathcal{F}, *)$-axial algebra.
\end{proof}

\section{Classifying fusion laws}
\label{sec:fusion_laws}

In order to classify graded $2$-generated axial algebras, we first need to classify the possible fusion laws that they can arise from.

\begin{prop}
    \label{prop:twoevalsfusionlaws}
    Suppose that $(\mathcal{F}, *)$ is a graded symmetric fusion law such that $\mathcal{F} = \{e, \alpha\}$ where $e$ is a unit. Then either $e*x = \emptyset$ for all $x \in \mathcal{F}$, or $(\mathcal{F}, *)$ is contained in the following fusion law.
    \begin{center}
    \def\arraystretch{1.5}
    \begin{tabular}{>{$}c<{$}|>{$}c<{$}>{$}c<{$}}
      & e & \alpha  \\ \hline
    e & e & \alpha \\
    \alpha & \alpha  & e \\
    \end{tabular}
    \end{center}
\end{prop}

\begin{proof}
    Let $\xi: \mathcal{F} \rightarrow G$ be a non-trivial adequate grading of $(\mathcal{F}, *)$. Suppose that there exists $x \in \mathcal{F}$ such that $e*x \neq \emptyset$. Then, as $e$ is a unit, $e*x = \{x\}$ and $\xi(e) = 1$. In addition, $e * e \subseteq  \{e\}$ and $e * \alpha \subseteq  \{\alpha\}$. It remains only to determine the value of $\xi(\alpha)$. If $\alpha \in \alpha*\alpha$ then $\xi(\alpha) \in \xi(\alpha)*\xi(\alpha) = \{\xi(\alpha)^2\}$ thus $ \xi(\alpha) = 1$, a contradiction as $\xi$ is not trivial. Therefore $(\mathcal{F}, *)$ must be contained in the fusion law given above.
\end{proof}

\begin{prop}
    \label{prop:threeevalsfusionlaws}
    Suppose that $(\mathcal{F}, *)$ is a graded symmetric fusion law such that $\mathcal{F} = \{e, \alpha, \beta\}$ where $e$ is a unit. Then either $e*x = \emptyset$ for all $x \in \mathcal{F}$, or $(\mathcal{F}, *)$ is contained in at least one of the following fusion laws.
    \begin{center}
        \begin{tabular}{llll}
            (a) & \def\arraystretch{1.5}
            \begin{tabular}{>{$}c<{$}|>{$}c<{$}>{$}c<{$}>{$}c<{$}}
                  & e & \alpha & \beta  \\ \hline
                e & e & \alpha & \beta \\
                \alpha & \alpha  & e, \alpha & \beta \\
                \beta & \beta  & \beta & e, \alpha \\
            \end{tabular} &
            (b) \def\arraystretch{1.5}
            \begin{tabular}{>{$}c<{$}|>{$}c<{$}>{$}c<{$}>{$}c<{$}}
                  & e & \alpha & \beta  \\ \hline
                e & e & \alpha & \beta \\
                \alpha & \alpha  & e, \beta & \alpha \\
                \beta & \beta  & \alpha & e, \beta \\
            \end{tabular} \\
            &&& \\
            (c) & \def\arraystretch{1.5}
            \begin{tabular}{>{$}c<{$}|>{$}c<{$}>{$}c<{$}>{$}c<{$}}
                  & e & \alpha & \beta  \\ \hline
                e & e & \alpha & \beta \\
                \alpha & \alpha  & e & e \\
                \beta & \beta  & e & e \\
            \end{tabular} &
            (d) \def\arraystretch{1.5}
            \begin{tabular}{>{$}c<{$}|>{$}c<{$}>{$}c<{$}>{$}c<{$}}
                  & e & \alpha & \beta  \\ \hline
                e & e & \alpha & \beta \\
                \alpha & \alpha  & \beta & e \\
                \beta & \beta  & e & \alpha \\
            \end{tabular}
        \end{tabular}
    \end{center}
    Note that fusion laws (a) and (b) are isomorphic.
\end{prop}

\begin{proof}
    Suppose that there exists $x \in \mathcal{F}$ such that $e*x \neq \emptyset$. Then, as $e$ is a unit, $e*x = \{x\}$ and $\xi(e) = 1$. In addition, $e * e \subseteq  \{e\}$ $e * \alpha \subseteq  \{\alpha\}$ and $e * \beta \subseteq  \{\beta\}$.

    We let $g := \xi(\alpha)$ and $h := \xi(\beta)$. Then
    \[
        \xi(\alpha * \alpha) \subseteq \xi(\alpha) * \xi(\alpha) = \{g^2\}.
    \]
    First assume that $\alpha \in \alpha * \alpha$. Then $g \in \xi(\alpha * \alpha) = \{g^2\}$ and so we must have $g = 1$. If also $\beta \in \alpha * \alpha$ then $\xi(e) = \xi(\alpha) = \xi(\beta) = 1$ and $\xi$ is trivial. Thus $\alpha * \alpha \subseteq \{e, \alpha\} $. As $\xi(\alpha) =1$, $\xi(\alpha *\beta) \subseteq \xi(\alpha)*\xi(\beta) = \{h\}$ and so $\alpha *\beta \subseteq \{\beta\}$. Finally, we cannot have $\beta \in \beta*\beta$ as otherwise $\xi$ would be trivial. Thus, if $\alpha \in \alpha * \alpha$, $(\mathcal{F}, *)$ is contained in case (a). Similarly, we can deduce that if $\beta \in \beta*\beta$ then $(\mathcal{F}, *)$ is contained in case (b). We can henceforth assume that $\alpha \not\in \alpha * \alpha$ and that $\beta \not\in \beta * \beta$.

    We now argue based on the possible values of $\alpha*\beta$. If $\alpha \in \alpha*\beta$ then $g = \xi(\alpha) \in \{gh\}$ and so $h = 1$. If also $e \in \alpha*\beta$ or $\beta \in \alpha*\beta$ then $g = h = 1$ and $\xi$ is trivial and so $\alpha *\beta \subseteq \{\alpha\}$. Finally, if $\alpha \in \beta *\beta$ then $g \in \xi(\beta)*\xi(\beta) = \{1\}$, again giving a contradiction. Thus, if $\alpha \in \alpha*\beta$, we are in case (b). Similarly, if $\beta \in \alpha*\beta$, we are in case (a). We can henceforth also assume that $\alpha \not\in \alpha * \beta$ and also that $\beta \not\in \alpha * \beta$. In particular, $\alpha * \beta \subseteq \{e\}$.

    If $\alpha * \beta = \{e\}$ then we must have $gh = 1$. Furthermore, if $\beta \in \alpha *\alpha$ then $h = g^2$ and so $g^3 = 1$ and if $e \in \alpha *\alpha$ then $g^2 = 1$. Thus if $\alpha *\alpha = \{e, \beta\}$ then $g = h = 1$, a contradiction as $\xi$ is non-trivial. Hence either $\alpha *\alpha \subseteq \{e\}$, or $\alpha *\alpha \subseteq \{\beta\}$. Similarly, we can also conclude that either $\beta *\beta \subseteq \{e\}$, or $\beta *\beta \subseteq \{\alpha\}$.

    If $\alpha * \alpha = \emptyset$ or $\beta * \beta = \emptyset$ then $(\mathcal{F}, *)$ is contained either in case (c) or (d), depending on the value of $\beta * \beta$ or $\alpha * \alpha$ respectively. Thus we can assume that $\alpha * \alpha \neq \emptyset$ and $\beta * \beta \neq \emptyset$.

    Now, if $\alpha * \alpha = \{e\}$ then $g^2 = 1$ and so, as $gh = 1$, $g = h$. If also $\beta *\beta = \{\alpha\}$ then $h^2 = g= h$ and so $g = h = 1$, a contradiction. Similarly, we cannot have $\beta * \beta = \{e\}$ and $\alpha * \alpha  = \{\beta\}$. If $\alpha * \alpha = \{e\}$ and $\beta * \beta = \{e\}$ then we are in case (c). Conversely, if $\alpha * \alpha = \{\beta\}$ and $\beta * \beta = \{\alpha\}$ then we are in case (d).

    Next we consider the case $\alpha * \beta = \emptyset$. We again look at possibilities for the set $\alpha *\alpha$. If $\alpha *\alpha = \{e, \beta\}$ then $1 = h = g^2$. If also $\alpha \in \beta *\beta$ then $g = h^2 = 1$ and so $\xi$ is trivial, a contradiction. Thus $\beta *\beta \subseteq \{e\}$ and we are in case (b). If $\alpha *\alpha = \{e\}$ then we are in case (a). If $\alpha *\alpha = \{\beta\}$ then we cannot have $\beta * \beta = \{e, \alpha\}$ as otherwise $ g= h= 1$, a contradiction. Thus either $\beta * \beta \subseteq \{e\}$ or $\beta * \beta \subseteq \{\alpha\}$  and we are in case (b) or (d) respectively. Finally, if $\alpha *\alpha = \emptyset$ then we are in case (b). This exhausts all possibilities for gradings of $(\mathcal{F}, *)$.
\end{proof}

Given a fusion law $(\mathcal{F}, *)$, recall from Proposition \ref{prop:mingrading} that the group $\Gamma_\mathcal{F}$ is defined as
\[
    \Gamma_{\mathcal{F}} := \langle \gamma_x : x \in \mathcal{F} \mid \gamma_x \gamma_y = \gamma_z \textrm{ whenever } z \in x *y \rangle.
\]
It is easy to calculate this group for the fusion laws classified in this section.

\begin{prop} \leavevmode
    \begin{enumerate}
        \item Suppose that $(\mathcal{F}, *)$ is the fusion law given in Proposition \ref{prop:twoevalsfusionlaws}. Then $\Gamma_{\mathcal{F}} \cong C_2$.
        \item Suppose that $(\mathcal{F}, *)$ is one of the fusion laws given in parts (a), (b) and (c) of Proposition \ref{prop:threeevalsfusionlaws}. Then $\Gamma_{\mathcal{F}} \cong C_2$.
        \item Suppose that $(\mathcal{F}, *)$ is the fusion laws given in part (d) of Proposition \ref{prop:threeevalsfusionlaws}. Then $\Gamma_{\mathcal{F}} \cong C_3$.
    \end{enumerate}
\end{prop}

More generally, the following result shows that sublaws of graded fusion laws are also graded.

\begin{prop}
    Suppose that $(\mathcal{F}, *_{\mathcal{F}})$ is a fusion law that admits a grading $f: (\mathcal{F}, *_{\mathcal{F}}) \rightarrow (G, *_G)$ for some group fusion law $(G, *_G)$ and suppose that $(\mathcal{H}, *_{\mathcal{H}})$ is a sublaw of $(\mathcal{F}, *_{\mathcal{F}})$. Then $f$ restricted to $H$ is also a $G$-grading for $(\mathcal{H}, *_{\mathcal{H}})$. Furthermore, if $f$ is adequate for $(\mathcal{F}, *_{\mathcal{F}})$ and $\mathcal{F} = \mathcal{H}$ then $f$ is also adequate for $(\mathcal{H}, *_{\mathcal{H}})$.
\end{prop}

\begin{proof}
    As $(\mathcal{H}, *_{\mathcal{H}})$ is a sublaw of $(\mathcal{F}, *_{\mathcal{F}})$, for all $x, y \in \mathcal{H}$
    \[
        f(x*_{\mathcal{H}}y) \subseteq f(x*_{\mathcal{F}}y) \subseteq f(x) *_G f(y)
    \]
    and so $f$ is also a $G$-grading of $(\mathcal{H}, *_{\mathcal{H}})$. If $f$ is adequate and $\mathcal{F} = \mathcal{H}$ then $f(\mathcal{F}) = f(\mathcal{H})$ generates $G$ and so $f$ is clearly also adequate for $(\mathcal{H}, *_{\mathcal{H}})$.
\end{proof}

\section{The two eigenvalue case}
\label{sec:2evals}

\begin{thm}
    \label{thm:2evals}
    Suppose that $U = \langle \langle a_0, a_1 \rangle \rangle$ is a primitive axial algebra with the following fusion law for some $\alpha \in \mathbb{C}$.
    \begin{center}
    \def\arraystretch{1.5}
    \begin{tabular}{>{$}c<{$}|>{$}c<{$}>{$}c<{$}}
      & 1 & \alpha  \\ \hline
    1 & 1 & \alpha \\
    \alpha & \alpha  & 1, \alpha \\
    \end{tabular}
    \end{center}
    Then either $U$ is the $1$-dimensional algebra $1A$, or $U$ is a quotient of the $2$-dimensional algebra $2B(\alpha)$.
\end{thm}

\begin{proof}
    Let $V$ be the universal primitive $2$-generated axial algebra with the fusion law as above. We first calculate that
    \begin{align*}
        v_1 &:= f_{\{\alpha\}}(ad_{a_0})(a_1) = a_0a_1 - \alpha a_1 \in V_1^{a_0} \\
        v_\alpha &:= f_{\{1\}}(ad_{a_0})(a_1) = a_0a_1 -  a_1 \in V_\alpha^{a_0}
    \end{align*}

    We denote the projections of $V$ onto the subspaces $\langle a_0 \rangle$ and $\langle a_1 \rangle$ as $\varphi_{a_0}$ and $\varphi_{a_1}$ respectively and let $x = \varphi_{a_0}(a_1)$ and $y = \varphi_{a_1}(a_0)$.

    As $a_1 = \frac{1}{1-\alpha}(v_1 - v_\alpha)$, from Proposition \ref{prop:projections}, we calculate that
    $\varphi_{a_0}(v_1) = (1 - \alpha)x$ and so, as $V$ is primitive,
    \[
        v_1 = a_0a_1 - \alpha a_1  = (1 - \alpha) x a_0.
    \]
    Thus
    \begin{equation}
        \label{eq:2evals}
        a_0a_1 =  (1 - \alpha) x a_0 + \alpha a_1.
    \end{equation}

    Let $\sigma$ be the semiautomorphism that exchanges $a_0$ and $a_1$ (see Proposition \ref{prop:sigma}). Then
    \begin{equation*}
        v := a_0a_1 - (a_0a_1)^\sigma = ((1-\alpha)x - \alpha)a_0 - ((1-\alpha)y - \alpha)a_1 = 0
    \end{equation*}
    As $v = 0$, $\varphi_{a_0}(v) = x - \alpha - (1 - \alpha)xy = 0$ and $\varphi_{a_1}(v) = (1 - \alpha)xy - y + \alpha = 0$. In particular, $\varphi_{a_0}(v) + \varphi_{a_1}(v) = x - y = 0$ so $x = y$. Then
    \[
        \varphi_{a_0}(v) = x - \alpha - (1 - \alpha)x^2 = ((\alpha - 1)x + \alpha)(x - 1) = 0.
    \]

    Now let $U$ be a primitive $2$-generated algebra with the same fusion law. As $U$ is a quotient of the universal algebra $V$, from $\varphi_{a_0}(v) = 0$ above, we have that in $U$, $x = \varphi_{a_0}(a_1) \in \{1, \frac{\alpha}{1 - \alpha}\}$.

    First suppose that $x = \frac{\alpha}{1 - \alpha}$. Then (\ref{eq:2evals}) becomes $a_0a_1 = \alpha(a_0 + a_1)$, and $U$ is the algebra $2B(\alpha)$. Next, if $x = 1$ then $v = (1 - 2 \alpha)(a_0 - a_1) = 0$. If $\alpha = \frac{1}{2}$ then (\ref{eq:2evals}) becomes $a_0a_1 = \frac{1}{2}(a_0 + a_1)$, and we are again in the $2B(\alpha)$ case. Otherwise, if $\alpha \neq \frac{1}{2}$ then we must have $a_0 = a_1$ and $U$ is a $1$-dimensional algebra.

%
\end{proof}

\begin{coro}
    \label{coro:2evalsgraded}
    Suppose that $U = \langle \langle a_0, a_1 \rangle \rangle$ is a primitive axial algebra with the following fusion law for some $\alpha \in \mathbb{C}$.
    \begin{center}
    \def\arraystretch{1.5}
    \begin{tabular}{>{$}c<{$}|>{$}c<{$}>{$}c<{$}}
      & 1 & \alpha  \\ \hline
    1 & 1 & \alpha \\
    \alpha & \alpha  & 1 \\
    \end{tabular}
    \end{center}
    Then $U$ is either the $1$-dimensional algebra $1A$ or is a quotient of the $2$-dimensional algebra $2B(\alpha)$. In the latter case, we must have $\alpha \in \{-1, \frac{1}{2}\}$.
\end{coro}

\begin{proof}
    This fusion law is a sublaw of that considered in Theorem \ref{thm:2evals}. Thus, from Theorem \ref{thm:sublaw}, either $U$ is the algebra $1A$ or is a quotient of the algebra $2B(\alpha)$.

    We can therefore assume that $V_\alpha^{a_0}$ is at most $1$-dimensional and is spanned by the vector
    \[
        v_\alpha = \alpha a_0 + (\alpha - 1) a_1.
    \]
    Using the products given by Theorem \ref{thm:2evals}, we calculate that
    \[
        v_\alpha v_\alpha = (2\alpha-1)(\alpha^2a_0 + (\alpha^2 - 1)a_1).
    \]

    In this fusion law, we must have $v_\alpha v_\alpha \in V_1^{a_0}$ and so $a_0(v_\alpha v_\alpha) - v_\alpha v_\alpha = 0$. We calculate that
    \[
        a_0(v_\alpha v_\alpha) - v_\alpha v_\alpha = (2\alpha - 1)(\alpha -1)(\alpha +1) (\alpha a_0 + (\alpha - 1)a_1).
    \]

    If $\alpha \in \{-1,  \frac{1}{2}\} $ then this vector is $0$ and we are done. Otherwise, as $\alpha \neq 1$, we can conclude that $a_1 = \frac{\alpha}{1-\alpha}a_0$ and so $U = \langle a_0 \rangle$ as required.
\end{proof}

The following result is not required for the main results; we include it only for completeness.

\begin{coro}
    \label{coro:jordansubtable}
    Suppose that $U = \langle \langle a_0, a_1 \rangle \rangle$ is a primitive axial algebra with the following fusion law for some $\alpha \in \mathbb{C}$.
    \begin{center}
    \def\arraystretch{1.5}
    \begin{tabular}{>{$}c<{$}|>{$}c<{$}>{$}c<{$}}
      & 1 & \alpha  \\ \hline
    1 & 1 & \alpha \\
    \alpha & \alpha  & \alpha \\
    \end{tabular}
    \end{center}
    Then $U$ is either the $1$-dimensional algebra $1A$ or is a quotient of the $2$-dimensional algebra $2B(\alpha)$. In the latter case, we must have $\alpha \in \{0, \frac{1}{2}\}$.
\end{coro}

\begin{proof}
    As in the proof of Corollary \ref{coro:2evalsgraded}, we can conclude that $v_\alpha = \alpha a_0 + (\alpha - 1) a_1$ and that
    \[
        v_\alpha v_\alpha = (2\alpha-1)(\alpha^2a_0 + (\alpha^2 - 1)a_1).
    \]
    In this fusion law, we must have $v_\alpha v_\alpha \in V_\alpha^{a_0}$ and so $a_0(v_\alpha v_\alpha) - \alpha v_\alpha v_\alpha = 0$. We calculate that
    \[
        a_0(v_\alpha v_\alpha) - v_\alpha v_\alpha = \alpha (2\alpha - 1)(\alpha -1) a_0.
    \]
    This must be equal to zero and so, as $\alpha \neq 1$ and $a_0$ is a non-zero vector, we conclude that $\alpha \in \{0, \frac{1}{2}\}$, as required.
\end{proof}

Finally, we classify all non-trivial proper ideals of the $2$-dimensional algebra given in Theorem \ref{thm:2evals}.

\begin{lem}
    \label{lem:ideal}
    Let $V$ be an algebra over a field $\mathbf{k}$ and let $a \in V$ be any non-zero vector. Additonally, suppose that there exists a basis $\{v_0, \dots, v_n\}$ of $V$ such that for $0 \leq i \leq n$, $a v_i = \alpha_i v_i$ where the $\alpha_i$ are pairwise distinct. Then if $I$ is a non-trivial ideal of $V$, we must have $v_i \in I$ for some $i$.
\end{lem}

\begin{proof}
    As $I$ is non-trivial, there exists $v \in I$ such that $v \neq 0$. If $v = v_i$ for some $0 \leq i \leq n$ then we are done. Otherwise we can assume that there exists non-zero scalars $\lambda_1, \dots, \lambda_m$ such that
    \[
        v = v_{i_0} + \sum_{j = 1}^m \lambda_j v_{i_j}
    \]
    where $0 \leq i_0 < i_1 < \dots < i_m \leq n$ and $m \geq 1$.

    Then also $av \in I$ so
    \[
        av - \alpha_{i_0} v = \sum_{j = 1}^m (\alpha_{i_j} - \alpha_{i_0})\lambda_j v_{i_j}  \in I.
    \]
    Then each of the scalars $(\alpha_{i_j} - \alpha_{i_0}) \lambda_j$ are non-zero, so we can proceed by induction to conclude that we must have $v_{i_n} \in I$, for example.
\end{proof}

\begin{prop}
    \label{prop:twoevalsideals}
    Suppose that $V = \langle \langle a_0, a_1 \rangle \rangle$ is the algebra $2B(\alpha)$ for some $\alpha \in \mathbb{C}$. If $\alpha \notin \{0, \frac{1}{2}\}$ then $V$ is simple. Otherwise, if $\alpha = \frac{1}{2}$ then $I = \langle a_0 - a_1 \rangle$ is an ideal of $V$ and if $\alpha = 0$ then both $I = \langle a_0 \rangle$ and $I = \langle a_1 \rangle$ are ideals. In general, if $I$ is a non-trivial proper ideal of $V$ then the algebra $V/I$ is the $1$-dimensional algebra $1A$.
\end{prop}

\begin{proof}
    Suppose that $V$ contains a proper, non-trivial ideal $I$ and suppose that $\alpha \neq 0$. Clearly $V$ obeys the conditions of Lemma \ref{lem:ideal} with basis $\{a_0, v_\alpha\}$ where
    \[
        v_\alpha = \alpha a_0 + ( 1- \alpha) a_1.
    \]

    First suppose that $a_0 \in I$. Then $a_0 v_\alpha = \alpha v_\alpha \in I$. Thus, as $I$ is proper, $a_0 \in I$ if and only if $\alpha = 0$. If $\alpha = 0$ then $\langle a_0 \rangle$ is an ideal.

    If $v_\alpha \in I$ then also $v_\alpha v_\alpha \in I$. By taking a linear combination of $v_\alpha$ and $v_\alpha v_\alpha$, we conclude that in this case
    \[
        (2\alpha - 1)\alpha a_0 \in I.
    \]
    As $I$ is proper, this implies that $\alpha \in \{0, \frac{1}{2}\}$. If $\alpha = 0$ then $\langle v_\alpha \rangle = \langle a_1 \rangle$ is an ideal and if $\alpha = \frac{1}{2}$ then $\langle v_\alpha \rangle = \langle a_0 - a_1 \rangle$ is an ideal. The final claim is easily verified for each value of $I$.
\end{proof}

\section{The three eigenvalue case}
\label{sec:3evals}

\begin{lem}
    \label{lem:prelims}
    Suppose that $V = \langle \langle a_0, a_1 \rangle \rangle$ is a primitive $(\mathcal{F}, *)$-axial algebra such that $\mathcal{F} = \{1, \alpha, \beta\}$ for some $\alpha, \beta \in \mathbb{C}$. Then for some $x,y \in \mathbb{C}$,
    \begin{align*}
        a_0(a_0a_1) &= (\alpha-1)(\beta-1)x a_0 - \alpha\beta a_1 + (\alpha + \beta) a_0a_1 \\
        a_1(a_0a_1) &= (\alpha-1)(\beta-1)y a_1 - \alpha\beta a_0 + (\alpha + \beta) a_0a_1.
    \end{align*}
    Moreover, there exist eigenvectors
    \begin{align*}
        v_\alpha &= (\beta-1)x a_0 - \beta a_1 + a_0a_1 \in V_\alpha^{a_0} \\
        v_\beta &= (\alpha-1)x a_0 - \alpha a_1 + a_0a_1 \in V_\beta^{a_0}.
    \end{align*}
    Finally,
    \begin{align*}
        \frac{v_\alpha v_\alpha - v_\beta v_\beta}{\alpha - \beta} =& \left((\alpha + \beta - 2\alpha\beta) x^2 - 2\alpha\beta\right) a_0 + \left(2\alpha \beta x + 2 (\alpha - 1)(\beta-1)y - \alpha - \beta \right)a_1 \\
        &- 2\left( x - \alpha - \beta \right)a_0a_1 \\
        \frac{v_\alpha v_\alpha - v_\alpha v_\beta}{\alpha -\beta} =& (\alpha(1-\beta)x^2-\alpha\beta)a_0 + (\alpha\beta x + (\alpha-1)(\beta-1)y - \beta)a_1 \\
        &+ ((\beta - \alpha - 1)x + \alpha + \beta)a_0a_1.
    \end{align*}
\end{lem}

\begin{proof}
    Let $a_1 = v'_1 + v'_\alpha + v'_\beta$. Then we have
    \begin{align*}
        v'_1 &= \frac{\alpha\beta a_1 - (\alpha + \beta)a_0a_1 + a_0(a_0a_1)}{(\alpha-1)(\beta-1)} \\
        v'_\alpha &= \frac{\beta a_1  - (\beta+1)a_0a_1 + a_0(a_0a_1)}{(\alpha-1)(\alpha -\beta)} \\
        v'_\beta &= \frac{\alpha a_1  - (\alpha+1)a_0a_1 + a_0(a_0a_1)}{(\beta-1)(\beta -\alpha)}
    \end{align*}

    As the algebra is primitive, from Proposition \ref{prop:projections} we must have $v'_1 = x a_0$ where $x = \varphi_{a_0}(a_1)$. Thus we have that
    \[
    a_0(a_0a_1) = (\alpha-1)(\beta-1)x a_0 - \alpha\beta a_1 + (\alpha + \beta) a_0a_1.
    \]
    Applying the semiautomorphism $\sigma$ that exchanges $a_0$ and $a_1$ gives
    \[
    a_1(a_0a_1) = (\alpha-1)(\beta-1)y a_1 - \alpha\beta a_0 + (\alpha + \beta) a_0a_1
    \]
    where $y = \varphi_{a_1}(a_0)$.

    Substituting this value into the $v'_\alpha$ and $v'_\beta$ above and rescaling gives the eigenvectors
    \begin{align*}
        v_\alpha &= (\beta-1)x a_0 - \beta a_1 + a_0a_1 \\
        v_\beta &= (\alpha-1)x a_0 - \alpha a_1 + a_0a_1
    \end{align*}

    Using these values we can calculate $v_\alpha v_\alpha - v_\beta v_\beta$ and $v_\alpha v_\alpha - v_\alpha v_\beta$.
\end{proof}

\begin{lem}
    \label{lem:projections}
    Let $x = \varphi_{a_0}(a_1)$ and $y = \varphi_{a_1}(a_0)$ as in Lemma \ref{lem:prelims}. Then
    \begin{enumerate}[i)]
        \item $\varphi_{a_0}(a_0a_1) = x$;
        \item $\varphi_{a_1}(a_0a_1) = y$;
        \item $\varphi_{a_1}(v_\alpha) = (\beta - 1)xy - \beta + y$;
        \item $\varphi_{a_1}(v_\beta) = (\alpha - 1)xy - \alpha + y$;
    \end{enumerate}
\end{lem}

\begin{proof}
    We first note that
    \begin{align*}
        \varphi_{a_0}(v_\beta) = \varphi_{a_0}(a_0a_1) - x = 0.
    \end{align*}
    So $\varphi_{a_0}(a_0a_1) = x$. By applying the semiautomorphism $\sigma$, we conclude that $\varphi_{a_1}(a_0a_1) = y$. We can then calculate $\varphi_{a_1}(v_\alpha)$ and $\varphi_{a_1}(v_\beta)$ as required.
\end{proof}

\subsection{Fusion law (a)}

\begin{thm}
    \label{thm:3evals}
    Suppose that $U = \langle \langle a_0, a_1 \rangle \rangle$ is a primitive axial algebra with the following fusion law for some $\alpha, \beta \in \mathbb{C}$.
    \begin{center}
        \def\arraystretch{1.5}
        \begin{tabular}{>{$}c<{$}|>{$}c<{$}>{$}c<{$}>{$}c<{$}}
              & 1 & \alpha & \beta  \\ \hline
            1 & 1 & \alpha & \beta \\
            \alpha & \alpha  & 1, \alpha & \beta \\
            \beta & \beta  & \beta & 1, \alpha \\
        \end{tabular}
    \end{center}
    Then $U$ is a quotient of one of the following.
    \begin{enumerate}[(i)]
        \item The $3$-dimensional algebra with basis $\{a_0, a_1, a_0a_1\}$ such that
        \begin{align*}
            (a_0a_1)(a_0a_1) =&  (\beta(\alpha-1)(\alpha + \beta - 1)x - \alpha\beta (\alpha + \beta))(a_0 + a_1) \\
            &+ ((1 - \alpha)(\alpha - \beta + 1)x +\alpha^2+3\alpha\beta + 2\beta^2-\beta)a_0a_1.
        \end{align*}
        where if $\beta \neq \frac{1}{2}$ then $x = \frac{\alpha+\beta}{2(1 -\alpha)}$ or, if $\beta = \frac{1}{2}$, then $x$ may take any value in $\mathbb{C}$.
        \item The $2$-dimensional algebra $2B(\alpha)$.
        \item The $1$-dimensional algebra $1A$.
    \end{enumerate}
\end{thm}

\begin{proof}
    Let $(\mathcal{F}, *)$ be the fusion law above and let $V$ be the universal $2$-generated primitive $(\mathcal{F}, *)$-axial algebra.

    In this case, $v_\alpha v_\alpha - v_\beta v_\beta \in V_1^{a_0} \oplus V_\alpha^{a_0}$. As $v_\alpha$ and $a_0$ are also in $V_1^{a_0} \oplus V_\alpha^{a_0}$, we can take a linear combination of these three vectors to conclude that
    \[
        p(x, y)a_1 \in V_1^{a_0} \oplus V_\alpha^{a_0}.
    \]
    where $x = \varphi_{a_0}(a_1)$ and $y = \varphi_{a_1}(a_0)$ and
    \[
        p(x, y) = 2\beta(\alpha-1)x + 2(\alpha-1)(\beta-1)y + (\alpha + \beta)(2\beta - 1).
    \]
    We can then check that
    \[
        x(\alpha - \beta)a_0 +  v_\alpha - (\alpha - \beta)a_1 = v_\beta.
    \]
    Thus, as $p(x,y)a_1 \in V_1^{a_0} \oplus V_\alpha^{a_0}$, $p(x,y)v_\beta \in (V_1^{a_0} \oplus V_\alpha^{a_0}) \cap V_\beta^{a_0} = 0$. In particular, $p(x,y)\varphi_{a_1}(v_\beta) = 0$ so either $p(x,y) = 0$ or $q(x,y) := \varphi_{a_1}(v_\beta) = 0$.

    We first assume that $p(x,y) = 0$. Then, by applying the semiautomorphism $\sigma$, we can also conclude that $p(y,x) = 0$ and, in particular,
    \[
        p(x, y) - p(y, x) = 2(x - y)(\alpha - 1)  = 0
    \]
    and so, as $\alpha \neq 1$, we can conclude that $x = y$. Otherwise, if $q(x,y) = 0$ then, using the value of $q(x,y) = \varphi_{a_1}(v_\beta)$ as given in Lemma \ref{lem:projections}, we similarly find that
    \[
        q(x,y) - q(y,x) = y - x = 0
    \]
    and we can again conclude that $x=y$.

    At this stage,
    \begin{equation}
        \label{eq:3evals}
        p(x,y)v_\beta = (2\beta - 1)(2(\alpha - 1)x + \alpha + \beta) v_\beta = 0.
    \end{equation}
    and
    \begin{equation}
        \label{eq:3evalspoly}
        p(x,y)q(x,y) = (2\beta - 1)(2(\alpha - 1)x + \alpha + \beta)((\alpha - 1) x + \alpha))(x-1) = 0
    \end{equation}

    Now note that
    \begin{align*}
        v_\alpha v_\alpha - v_\alpha v_\beta = &\alpha(\alpha - \beta)((\beta - 1)x + \beta)(x -1) a_0 + (\alpha + \beta -1 )((1-\beta)x - \beta) v_\beta \\
         &+ ((3\alpha\beta -\alpha^2-2\alpha-\beta+1)x  +\alpha^2 + \alpha\beta - \beta) v_\alpha.
    \end{align*}
    As $a_0, v_\alpha, v_\alpha v_\alpha \in V_1^{a_0} \oplus V_\alpha^{a_0}$ and $v_\alpha v_\beta, v_\beta \in V_\beta^{a_0}$,
    \[
        v_\alpha v_\beta - (\alpha + \beta -1 )((1-\beta)x - \beta) v_\beta \in (V_1^{a_0} \oplus V_\alpha^{a_0})  \cap V_\beta^{a_0}
    \]
    and so $v_\alpha v_\beta - (\alpha + \beta -1 )((1-\beta)x - \beta) v_\beta = 0$. From this, we can calculate that
    \begin{align*}
        (a_0a_1)(a_0a_1) =&  (\beta(\alpha-1)(\alpha + \beta - 1)x - \alpha\beta (\alpha + \beta))(a_0 + a_1) \\
        &+ ((1 - \alpha)(\alpha - \beta + 1)x +\alpha^2+3\alpha\beta + 2\beta^2-\beta)a_0a_1.
    \end{align*}

    Now let $U$ be a $2$-generated primitive $(\mathcal{F}, *)$-axial algebra. From (\ref{eq:3evalspoly}), either $\beta = \frac{1}{2}$ or, in the algebra $U$, we have $x \in \{\frac{\alpha + \beta}{2(1 - \alpha)}, \frac{\alpha}{1 - \alpha}, 1\}$.

    If $\beta = \frac{1}{2}$ or if $\beta \neq \frac{1}{2}$ and $x = \frac{\alpha + \beta}{2(1 - \alpha)}$ then we are in case (i). We can check that this algebra is indeed a primitive $(\mathcal{F}, *)$-axial algebra as required. We can thus assume that $\beta \neq \frac{1}{2}$ and $x \neq \frac{\alpha + \beta}{2(1 - \alpha)}$. If $x = \frac{\alpha}{1 - \alpha}$ then (\ref{eq:3evals}) becomes $(2\beta - 1)(\beta - \alpha)v_\beta = 0$. Thus, as $\beta \neq \frac{1}{2}$ and $\alpha \neq \beta$, we have $v_\beta = 0$ and so $a_0a_1 = \alpha(a_0 + a_1)$ and we are in case (ii).

    If $x = 1$ then (\ref{eq:3evals}) becomes $(2\beta - 1)(3\alpha + \beta - 2)v_\beta = 0$. If $\beta \neq 2 - 3\alpha$ then again $v_\beta = 0$ and in this case we have
    \[
        a_0a_1 = (1 - \alpha)a_0 + \alpha a_1.
    \]

    We can now calculate the adjoint action of the axis $a_1$ and see that it has eigenvalues $1$ and $1-\alpha$ with corresponding eigenvectors $a_1$ and $a_0 - a_1$ respectively. As $a_1$ must be an axis with the above fusion law, we must have either $1 - \alpha = \alpha$ or $1 - \alpha = \beta$. If $1 - \alpha = \alpha$ then $\alpha = \frac{1}{2}$ and $a_0a_1 = \alpha(a_0 + a_1)$ so we are in case (ii). Thus we can assume that $\alpha \neq \frac{1}{2}$. The vector $a_0 - a_1 \in V_\beta^{a_1}$ but
    \[
        (a_0 - a_1)(a_0 - a_1) = (1 - 2\alpha)(a_0 - a_1) \in (V_1^{a_1} \oplus V_\alpha^{a_1}) \cap V_\beta^{a_1} = 0.
    \]
    So, as $\alpha \neq \frac{1}{2}$, $a_0 = a_1$ and we are in case (iii).

    The final case to consider is that of $x = 1$ and $\beta = 2 - 3\alpha$. In this case, as $\beta \neq \alpha$, $\beta \neq \frac{1}{2}$ and moreover,
    \[
        \frac{\alpha + \beta}{2(1 - \alpha)} = 1
    \]
    so in fact we are in case (i).
\end{proof}

\begin{prop}
    \label{prop:threeevalsideals}
    Let $V$ be the $3$-dimensional algebra defined in part (i) of Theorem \ref{thm:3evals}. Then all non-trivial proper ideals of $V$ and the values of $\alpha$, $\beta$ and $x$ for which they occur are given by Table \ref{tab:threeevalsideals}. If $I$ is a non-trivial proper ideal of $V$ then the algebra $V/I$ is either the $1$-dimensional algebra $1A$, or is one of the $2$-dimensional algebras $2B(\alpha)$ or $2B(\beta)$.
\end{prop}

\begin{table}
\begin{center}
\vspace{0.35cm}
\noindent
\begin{tabular}{>{$}c<{$}|>{$}l<{$}}
 I & \textrm{Parameters for which } I \textrm{ is an ideal}     \\
\hline
\langle v_\alpha \rangle & \beta = -1 \textrm{ and } x = - \frac{1}{2} \\
\langle v_\alpha \rangle & \beta = \frac{1}{2} \textrm{ and } x = 1  \\
& \\
\langle v_\beta \rangle & \beta = \frac{1}{2} \textrm{ and } x = \frac{\alpha}{\alpha - 1}  \\
& \\
\langle a_0, v_\alpha \rangle & \beta = 0   \\
& \\
\langle a_0, v_\beta \rangle & \alpha = 0 \textrm{ and } \beta = \frac{1}{2} \textrm{ and } x = 0   \\
& \\
\langle v_\alpha, v_\beta \rangle & \beta = \frac{1}{2}  \textrm{ and } x = 1   \\
\langle v_\alpha, v_\beta \rangle & \beta = 2 - 3\alpha \textrm{ and } x = 1   \\
\langle v_\alpha, v_\beta \rangle & \beta = 0 \textrm{ and } x = \frac{\alpha}{1 - \alpha}   \\
\end{tabular}
\caption{The non-trivial proper ideals of $V$}
\label{tab:threeevalsideals}
\end{center}
\end{table}

\begin{proof}
    Let $I$ be a non-trivial, proper ideal of $V$. The basis $\{a_0, v_\alpha, v_\beta\}$ of $V$ clearly satisfies the conditions of Lemma \ref{lem:ideal}.  Thus we can assume that at least one of $\{a_0, v_\alpha, v_\beta\}$ is contained inside $I$.

    If we were to have $a_0 \in I$ then also $a_0 v_\alpha = \alpha v_\alpha \in I$ and $a_0 v_\beta = \beta v_\beta \in I$ thus, as $I$ is a proper ideal, if $\alpha \neq 0$ and $\beta \neq 0$, $a_0 \notin I$.

    Now suppose that $a_0 \in I$ and $\alpha = 0$ then $\beta \neq 0$ and so we can conclude that $v_\beta \in I$. Then also $v_\beta v_\beta \in I$. By taking a linear combination of $a_0$, $v_\beta$ and $v_\beta v_\beta$, we conclude that $\beta ( 1- \beta) x a_1 \in I$. If $a_0, a_1 \in I$ then we must have $I = V$, a contradiction. Thus, as $\beta \notin \{1, 0\}$, $x = 0$. In this case, $\langle a_0, v_\beta \rangle$ is an ideal. Next, suppose that $a_0 \in I$ and $\beta = 0$ then $\alpha \neq 0$ and so we can conclude that $v_\alpha \in I$. From the fusion law, we see that $v_\alpha v_\alpha \in V_1^{a_0} \oplus V_\alpha^{a_0} = \langle a_0, v_\alpha \rangle$ for all $v \in V$. Thus $\langle a_0, v_\alpha \rangle$ is an ideal.

    If $v_\beta \in I$ then also $v_\beta v_\beta \in I$. We first suppose that $\beta = \frac{1}{2}$ and calculate that
    \[
        v_\beta v_\beta = \frac{1}{4}(2\alpha-1)((\alpha - 1)x + \alpha)(a_0 + a_1 - 2a_0a_1).
    \]
    As $\beta = \frac{1}{2}$, $\alpha \neq \frac{1}{2}$ and so $v_\beta v_\beta = 0$ if and only if $ x = \frac{\alpha}{\alpha - 1}$. In this case, $\langle v_\beta \rangle$ is an ideal. If additionally $v_\alpha \in I$ then also $v_\alpha v_\alpha \in I$. By taking a linear combination of $v_\alpha$ and $v_\alpha v_\alpha$, we conclude that also
    \[
        \alpha (2 \alpha - 1)^2 a_0 \in I.
    \]
    As $I$ is proper and contains $v_\alpha$ and $v_\beta$, we must have $a_0 \notin I$. As $\beta = \frac{1}{2}$, $\alpha \neq \frac{1}{2}$ and we can conclude that in this case $\alpha = 0$ (and therefore $x=0$) and $\langle v_\alpha, v_\beta \rangle$ is an ideal.

    Finally, if $v_\beta v_\beta \neq 0$ then $v := a_0 + a_1 - 2a_0a_1 \in I$. Then
    \[
        a_0 v - \alpha v = \frac{1}{2}(1 - \alpha)(x - 1)a_0 \in I.
    \]
    So if $x = 1$ then $v \in V_\alpha^{a_0}$ and we can check that $\langle v_\beta, v_\beta v_\beta \rangle = \langle v_\beta, v_\alpha \rangle$ is an ideal. Otherwise, if $x \neq 1$ then $v_\beta, a_0 \in I $ so we must be in the case $\alpha = 0$, $x = 0$ and $I = \langle v_\beta, a_0 \rangle$.

    We now suppose that $x = \frac{\alpha+\beta}{2(1 - \alpha)}$ and calculate that
    \[
        v_\beta v_\beta = (\alpha - \beta)((4\alpha \beta - \alpha - \beta)a_0 + 2(\alpha + \beta - 1)(\beta a_1 - a_0a_1)).
    \]
    In this case, we see that $v_\beta v_\beta = 0$ if and only if $\alpha = \beta = \frac{1}{2}$, so this cannot occur and we have $v_\beta v_\beta \neq 0$.

    We further calculate that
    \[
        a_0 (v_\beta v_\beta) - \alpha v_\beta v_\beta = (\alpha - \beta)^2 \beta (2 - 3 \alpha - \beta)a_0 \in I.
    \]
    So if $\beta = 0$ or $\beta = 2 - 3\alpha$ then $v_\beta v_\beta \in V_\alpha^{a_0}$ and we can check that $\langle v_\beta, v_\beta v_\beta \rangle = \langle v_\beta, v_\alpha \rangle$ is an ideal. Otherwise, again we must be in the case $\alpha = 0$, $x = 0$. However, in this case, this would imply that $\beta = 0$, a contradiction.

    If $v_\alpha \in I$ then also $v_\alpha v_\alpha \in I$. Then by taking a linear combination of $v_\alpha$ and $v_\alpha v_\alpha$, we obtain that
    \[
        \alpha (\alpha - \beta) ((\beta - 1)x + \beta)(x-1)a_0 \in I.
    \]

    If $\alpha = 0$ or $x \in \{1, \frac{\beta}{1 - \beta}\}$ then $v_\alpha v_\alpha \in V_\alpha^{a_0}$ and so $I = \langle v_\alpha \rangle$ is an ideal. Otherwise, $a_0 \in I$ and so we must be in the case where $\beta = 0$ and $I = \langle a_0, v_\alpha \rangle$.

    We check that in each case the algebra $V/I$ is either the $1$-dimensional algebra $1A$ or is one of the $2$-dimensional algebras $2B(\alpha)$ or $2B(\beta)$ (if $I = \langle v_\beta \rangle$ or  $I = \langle v_\alpha \rangle$ respectively).
\end{proof}

\subsection{Fusion law (c)}

\begin{thm}
    \label{thm:3evalsother}
    Suppose that $U= \langle \langle a_0, a_1 \rangle \rangle$ is a primitive axial algebra with the following fusion law for some $\alpha, \beta \in \mathbb{C}$.
    \begin{center}
        \def\arraystretch{1.5}
        \begin{tabular}{>{$}c<{$}|>{$}c<{$}>{$}c<{$}>{$}c<{$}}
              & 1 & \alpha & \beta  \\ \hline
            1 & 1 & \alpha & \beta \\
            \alpha & \alpha  & 1 & 1 \\
            \beta & \beta  & 1 & 1 \\
        \end{tabular}
    \end{center}
    Then $U$ is either the $1$-dimensional algebra $1A$ or is a quotient either of $2B(\alpha)$ or $2B(\beta)$. In the latter case then either $\alpha \in \{-1, \frac{1}{2}\}$ or $\beta \in \{-1, \frac{1}{2}\}$, respectively.
\end{thm}

\begin{proof}
    Let $(\mathcal{F}, *)$ be the fusion law above and let $V$ be the universal $2$-generated primitive $(\mathcal{F}, *)$-axial algebra.

    First, using Lemma \ref{lem:projections}, we calculate that $\frac{1}{\alpha - \beta}\varphi_{a_0}(v_\alpha v_\alpha - v_\alpha v_\beta) = p(x,y)$ where
    \[
        p(x,y) = (\beta - 1) x^2 + (\alpha-1)(\beta-1)xy + \alpha x - \alpha\beta.
    \]
    Then, as $V$ is primitive, $v := \frac{1}{\alpha - \beta}(v_\alpha v_\alpha - v_\alpha v_\beta) - p(x,y)a_0 = 0$. We can also calculate that
    \[
        \varphi_{a_1}(v) = (\alpha+1)(1 - \beta)x^2y + (\beta - 1)(1-\alpha)xy^2- (2\alpha - \beta + 1)xy + \alpha\beta x + (\alpha\beta + 1)y - \beta = 0.
    \]

    We now notice that exchanging $\alpha$ and $\beta$ preserves the fusion law and so, from Proposition \ref{prop:isomorphicfusionlaws}, induces an automorphism $\tau$ on $V$. Thus we can conclude that
    \[
        \varphi_{a_1}(v) - \varphi_{a_1}(v)^\tau = (\alpha - \beta)(2x^2y-3xy + 1) = 0
    \]
    At this stage, we apply the semiautomorphism $\sigma$ to conclude that
    \[
        \varphi_{a_1}(v) - \varphi_{a_1}(v)^\tau - [\varphi_{a_1}(v) - \varphi_{a_1}(v)^\tau]^\sigma = 2(\alpha - \beta)(x-y)xy = 0.
    \]
    As $\alpha \neq \beta$, either $x = y$ or $xy = 0$. If $xy = 0$ then $\varphi_{a_1}(v) - \varphi_{a_1}(v)^\tau = (\alpha - \beta) = 0$, a contradiction, and so $x = y$. We now have
    \begin{equation}
        \label{eq:3evalsother1}
        \varphi_{a_1}(v) - \varphi_{a_1}(v)^\tau = (\alpha - \beta)(x-1)^2(2x-1) = 0
    \end{equation}
    and
    \begin{equation}
        \label{eq:3evalsother2}
        \varphi_{a_1}(v) = (2\alpha x - 1)((1- \beta)x - \beta)(x - 1) = 0.
    \end{equation}

    Now let $z := \varphi_{a_0}(v_\alpha v_\alpha)$. Then, as $v_\alpha v_\alpha \in V_1^{a_0}$ and $V$ is primitive, $v_\alpha v_\alpha = z a_0$. We can use this to calculate that
    \[
        (a_0a_1)(a_0a_1) = (z + (1- 2\alpha)(\beta-1)^2x - 2\alpha\beta^2)a_0 + (2\beta(2\alpha-1)(\beta-1)x - \beta^2)a_1 + 2(\alpha(1 - \beta)x + \beta(\alpha +1))a_0a_1.
    \]

    We can now calculate the value of $v_\alpha v_\beta$, $v_\beta v_\beta$ and also of   $\varphi_{a_0}(v_\alpha v_\beta)$ and $\varphi_{a_0}(v_\beta v_\beta)$. Then
    \begin{align}
        \label{eq:vavb}
        \frac{v_\alpha v_\beta - \varphi_{a_0}(v_\alpha v_\beta) a_0}{\alpha - \beta} =& \alpha((2 \beta - 1)x + 1)xa_0 - ((2\alpha\beta - \alpha - \beta + 1)x - \beta)a_1 \nonumber \\
        & + ((\alpha - \beta + 1)x - \alpha - \beta)a_0a_1 = 0.
    \end{align}
    and
    \begin{align}
        \label{eq:vbvb}
        \frac{v_\beta v_\beta - \varphi_{a_0}(v_\beta v_\beta) a_0}{\alpha - \beta} =&
        (2(2\alpha\beta - \alpha - \beta)x + \alpha + \beta)xa_0 + (\alpha + \beta -  2(2\alpha\beta -\alpha - \beta - 1)x)a_1 \nonumber \\
        &+ 2(x - \alpha - \beta)a_0a_1 = 0.
    \end{align}

    Now let $U$ be a $2$-generated primitive $(\mathcal{F}, *)$-axial algebra. From (\ref{eq:3evalsother1}), in $U$, $x \in \{1, -\frac{1}{2}\}$.

    We first suppose that $x = -\frac{1}{2}$. In this case (\ref{eq:3evalsother2}) gives
    \[
        \varphi_{a_0}(v) = -\frac{3}{4}(\alpha + 1)(\beta + 1) = 0
    \]
    so either $\alpha = -1$ or $\beta = -1$. Without loss of generality, we assume that $\beta = -1$. Then (\ref{eq:vavb}) becomes
    \[
        -\frac{3\alpha}{2}(a_0 + a_1 + a_0a_1) = 0
    \]
    and (\ref{eq:vbvb})
    \[
        (1 - 2\alpha)(a_0 + a_1 + a_0a_1) = 0.
    \]
    Taking a linear combination of these two equations, we can conclude that $a_0a_1 = -(a_0 + a_1)$ and $U$ is the algebra $2B(\beta)$.

    We now suppose that $x = 1$. Then (\ref{eq:vavb}) becomes
    \[
        (2\beta - 1)(\alpha a_0 + (1- \alpha)a_1 - a_0a_1) = 0.
    \]
    We first suppose that $\beta \neq \frac{1}{2}$ then $a_0a_1 = \alpha a_0 + (1- \alpha)a_1$. We can now calculate the adjoint action of the axis $a_1$ and see that it has eigenvalues $1$ and $1-\alpha$ with corresponding eigenvectors $a_1$ and $a_0 - a_1$ respectively. As $a_1$ must be an axis with the above fusion law, we must have either $1 - \alpha = \alpha$ or $1 - \alpha = \beta$.

    If $1 - \alpha = \alpha$ then $\alpha = \frac{1}{2}$ and $a_0a_1 = \alpha(a_0 + a_1)$ so $U$ is the algebra $2B(\alpha)$. Thus we can assume that $\alpha \neq \frac{1}{2}$. The vector $a_0 - a_1 \in V_\beta^{a_1}$ but
    \[
        (a_0 - a_1)(a_0 - a_1) = (1 - 2\alpha)(a_0 - a_1) \in V_1^{a_1} \cap V_\beta^{a_1} = 0.
    \]
    So, as $\alpha \neq \frac{1}{2}$, $a_0 = a_1$ and $U$ is the algebra $1A$.

    Otherwise, if $\beta = \frac{1}{2}$ then (\ref{eq:vbvb}) becomes
    \[
        (\alpha - \frac{1}{2})(a_0 + a_0 - 2a_0a_1) = 0.
    \]
    As $\beta = \frac{1}{2}$, $\alpha \neq \frac{1}{2}$ and so $a_0a_1 = \frac{1}{2}(a_0 + a_1)$ and $U$ is the algebra $2B(\beta)$.
\end{proof}

\subsection{Fusion law (d)}

\begin{thm}
    \label{thm:C3graded}
    Suppose that $U = \langle \langle a_0, a_1 \rangle \rangle$ is a primitive axial algebra with the following fusion law for some $\alpha, \beta \in \mathbb{C}$.
    \begin{center}
        \def\arraystretch{1.5}
        \begin{tabular}{>{$}c<{$}|>{$}c<{$}>{$}c<{$}>{$}c<{$}}
              & 1 & \alpha & \beta  \\ \hline
            1 & 1 & \alpha & \beta \\
            \alpha & \alpha  & \beta & 1 \\
            \beta & \beta  & 1 & \alpha \\
        \end{tabular}
    \end{center}

    If $\alpha = \frac{1}{2}$ then either $U$ is the $1$-dimensional algebra $1A$ or $U$ is a quotient of the $3$-dimensional algebra with basis $\{a_0, a_1, a_0a_1\}$ such that
    \[
        (a_0a_1)(a_0a_1) = (\frac{1}{4} - \beta)(a_0 + a_1) + (2\beta +  \frac{1}{2})a_0a_1.
    \]

    If $\beta = \frac{1}{2}$ then either $U$ is the $1$-dimensional algebra $1A$ or $U$ is a quotient of the $3$-dimensional algebra with basis $\{a_0, a_1, a_0a_1\}$ such that
    \[
        (a_0a_1)(a_0a_1) = (\frac{1}{4} - \alpha)(a_0 + a_1) + (2\alpha +  \frac{1}{2})a_0a_1.
    \]

    Otherwise, if $\frac{1}{2} \notin \{\alpha, \beta\}$ then $U$ is the $1$-dimensional algebra $1A$.
\end{thm}

\begin{lem}
    \label{lem:C3graded}
    Let $V$ be the universal $2$-generated primitive axial algebra. Then in the coefficient ring of $V$:
    \begin{enumerate}[i)]
        \item $ x= y$;
        \item if $\beta \neq -\alpha -1$ then $x = y = 1$;
        \item if $\beta = -\alpha -1$ then $(2(2\alpha^2+2\alpha-1)x-1)((\alpha+2)x + \alpha + 1)(x-1) = 0$.
    \end{enumerate}
\end{lem}

\begin{proof}
    In this case, $v_\alpha v_\alpha - v_\alpha v_\beta \in V_1^{a_0} \oplus V_\beta^{a_0}$. As $v_\beta$ and $a_0$ are also in $V_1^{a_0} \oplus V_\beta^{a_0}$, we can take a linear combination of these vectors to conclude that
    \[
        p(x,y) a_1 \in V_1^{a_0} \oplus V_\beta^{a_0}
    \]
    where
    \[
        p(x,y) = (2\alpha\beta - \beta^2 - \beta)x + (\alpha -1)(\beta - 1)y + \beta^2 + \alpha \beta -\alpha.
    \]

    We can then check that
    \[
        x(\beta - \alpha)a_0 +  v_\beta + (\alpha - \beta)a_1 = v_\alpha.
    \]
    Thus, as $p(x,y)a_1 \in V_1^{a_0} \oplus V_\beta^{a_0}$, $p(x,y)v_\alpha \in (V_1^{a_0} \oplus V_\beta^{a_0}) \cap V_\alpha^{a_0} = 0$. In particular, $p(x,y)\varphi_{a_1}(v_\alpha) = 0$ and so either $p(x,y) = 0$ or $q(x,y) := \varphi_{a_1}(v_\alpha) = 0$.

    We first assume that $p(x,y) = 0$. We note that exchanging $\alpha$ and $\beta$ preserves the fusion law and so, from Proposition \ref{prop:isomorphicfusionlaws}, induces an automorphism $\tau$ of $V$. Thus
    \[
        p(x,y) - p(x,y)^\tau = (\alpha - \beta)(\alpha + \beta + 1)(x - 1) = 0.
    \]
    As $\alpha \neq \beta$, either $x = 1$ or $\beta = - \alpha - 1$. If $x = 1$ then by applying the semiautomorphism $\sigma$, we can also conclude that $y = 1$.

    If $\beta = - \alpha - 1$ then
    \[
        p(x,y) = 3\alpha(\alpha + 1)x + (\alpha + 2)(\alpha - 1)y - 1 = 0
    \]
    and also
    \[
        p(x,y) - p(x,y)^\sigma = 2(\alpha^2 + \alpha + 1)(x-y) =0.
    \]
    We could have equivalently considered $\alpha = - \beta - 1$ and concluded that
    \[
        p(x,y) - p(x,y)^\sigma = 2(\beta^2 + \beta + 1)(x-y) =0.
    \]
    Thus, either $x = y$ or $\alpha^2 + \alpha + 1 = 0$ and $\beta^2 + \beta + 1 = 0$.

    The roots of the polynomial $z^2 + z + 1$ are the complex numbers $\{\frac{1 \pm i\sqrt{3}}{2} \}$. Thus, if $x \neq y$ then we must have $\alpha, \beta \in \{\frac{1 \pm i\sqrt{3}}{2} \}$. However, as $\alpha \neq \beta$ we can then conclude that
    \[
        \alpha + \beta = \frac{1 + i\sqrt{3}}{2} + \frac{1 - i\sqrt{3}}{2} = 1.
    \]
    This is in contradiction with the fact that $\alpha + \beta = -1$ and so this case cannot occur. Thus, if $\beta = -\alpha-1 $ then we can also conclude that $x = y$.

    Now assume that $q(x,y) = 0$ where the value of $q(x,y) = \varphi_{a_1}(v_\beta)$ is given as in Lemma \ref{lem:projections}. Let $\sigma$ be the semiautomorphism of $V$ that exchanges $a_0$ and $a_1$. Then as $q(x,y) = 0$, also $q(x,y)^\sigma = 0$ and, in particular,
    \[
        q(x,y) - q(x,y)^\sigma = y - x =0.
    \]
    So we can conclude that $x = y$. Then
    \[
        q(x,y) = (\beta - 1)x^2 - \beta + x = ((\beta - 1)x + \beta)(x - 1) = 0
    \]
    so either $x = y = 1$ or $x = y = \frac{\beta}{1-\beta}$.

    Thus, at this stage, if $\beta \neq -\alpha-1$ then either $x = y = 1$, or $x = y = \frac{\beta}{1 - \beta}$. By applying the automorphism $\tau$ that exchanges $\alpha$ and $\beta$ then we can also conclude that either $x = y = 1$, or $x = y = \frac{\alpha}{1 - \alpha}$. As $\alpha \neq \beta$, we cannot have that $ x= y = \frac{\alpha}{1 - \alpha}$ and $x = y = \frac{\beta}{1 - \beta}$. Thus $x = y = 1$.

    Finally, if $\beta = -\alpha -1$ then
    \[
        p(x,y)q(x,y) = (2\alpha+1)(2(2\alpha^2+2\alpha-1)x-1)((\alpha+2)x + \alpha + 1)(x-1) = 0.
    \]
    As $\beta = -\alpha - 1$ and $\alpha \neq \beta$, we must have $\alpha \neq -\frac{1}{2}$ and so part (iii) holds.
\end{proof}

\begin{lem}
    \label{lem:C3gradedtrivial}
    Suppose that $V$ is as above and that $\beta \neq - \alpha - 1$ and $\frac{1}{2} \notin \{\alpha, \beta\}$. Then $V$ is 1-dimensional.
\end{lem}

\begin{proof}
    As in the proof of Lemma \ref{lem:C3graded}, we can show that
    \[
        p(x,y)v_\alpha \in (V_1^{a_0} \oplus V_\beta^{a_0}) \cap V_\alpha^{a_0} = 0.
    \]
    From Lemma \ref{lem:C3graded}, as $\beta \neq -\alpha - 1$, $x = y = 1$ and so
    \[
        p(x,y) = (2\alpha - 1)(2\beta -1).
    \]

    As $\frac{1}{2} \notin \{\alpha, \beta\}$ then $p(x,y) \neq 0$ and so we must conclude that $v_\alpha = 0$. In particular, this implies that
    \[
        a_0a_1 = (1 - \beta)a_0 + \beta a_1.
    \]
    Then, using the semiautomorphism $\sigma$ that exchanges $a_0$ and $a_1$, we conclude that
    \[
        a_0a_1 - (a_0a_1)^\sigma = (1 - 2\beta)(a_0 - a_1) = 0.
    \]
    As $\beta \neq \frac{1}{2}$, $a_0 = a_1$ and $V = \langle a_0 \rangle$ is $1$-dimensional.
\end{proof}

\begin{proof}[Proof of Theorem \ref{thm:C3graded}]
    From Lemmas \ref{lem:C3graded} and \ref{lem:C3gradedtrivial}, we have two cases to consider:
    \begin{enumerate}[i)]
        \item $\beta = -\alpha - 1$ and $x = y$;
        \item $\beta \neq -\alpha - 1$, $x = y = 1$ and $\frac{1}{2} \in \{\alpha, \beta\}$.
    \end{enumerate}
    In the latter case, we can assume without loss of generality that $\beta = \frac{1}{2}$, in which case $\alpha \neq -\frac{3}{2}$.

    We first note that
    \begin{align*}
        \label{eq:vavavbvb}
        v_\alpha v_\alpha - v_\beta v_\beta =& (\alpha - \beta)((2\alpha\beta - \alpha - \beta)x +2\alpha\beta)(x-1)a_0 \\
        &+ (2\alpha -1)(2(\beta - 1)x + \alpha + \beta ) v_\alpha -  (2\beta -1)(2(\alpha - 1)x + \alpha + \beta ) v_\beta.
    \end{align*}

    As $a_0, v_\alpha, v_\beta v_\beta \in V_1^{a_0} \oplus V_\alpha^{a_0}$ and $v_\beta, v_\alpha v_\alpha \in V_\beta^{a_0}$,
    \[
        v_\alpha v_\alpha + (2\beta -1)(2(\alpha - 1)x + \alpha + \beta ) v_\beta \in (V_1^{a_0} \oplus V_\alpha^{a_0}) \cap V_\beta^{a_0}
    \]
    and so $v_\alpha v_\alpha + (2\beta -1)(2(\alpha - 1)x + \alpha + \beta ) v_\beta = 0$. From this we can calculate the value of $(a_0a_1)(a_0a_1)$.

    If $x = y = 1$ and $\beta = \frac{1}{2}$ then
    \[
        (a_0a_1)(a_0a_1) = (\frac{1}{4} - \alpha)(a_0 + a_1) + (2\alpha +  \frac{1}{2})a_0a_1.
    \]
    We can check that these values give an axial algebra, as required.

    If $\beta = - \alpha - 1$ then
    \begin{align}
        \label{eq:C3graded}
        (a_0a_1)(a_0a_1) =& ((2\alpha^3 - 9\alpha^2 - 12\alpha+ 10)x^2 - (\alpha-1)(2\alpha+3)x -  2\alpha(\alpha+1)^2)a_0 \nonumber \\&
        + (4(2\alpha^2 + 2\alpha - 1)(x+1)  - \alpha(\alpha+1))a_1 + (6(\alpha^2 + \alpha - 1)x -1)a_0a_1.
    \end{align}

    Using this, we can calculate that $a_0(v_\beta v_\beta) - \alpha v_\beta v_\beta = r(x,y)a_0$ where
    \[
        r(x,y) = ((2\alpha^2+2\alpha-1)x + 2\alpha^2 + 2\alpha)(x-1).
    \]
    From the fusion law, we can conclude that $r(x,y) = 0$ or equivalently, that
    \[
         (2\alpha^2+2\alpha-1)x(x-1) = -(2\alpha^2 + 2\alpha)(x-1)
    \]
    By substituting this expression into part (iii) of Lemma \ref{lem:C3graded} we obtain that
    \[
        (2\alpha + 1)^2 ((\alpha+2)x + \alpha + 1)(x-1) = 0.
    \]
    As $\beta = -\alpha - 1$ and $\beta \neq \alpha$, $\alpha \neq - \frac{1}{2}$ and so we can conclude that
    \[
        ((\alpha+2)x + \alpha + 1)(x-1) = 0.
    \]
    Similarly, as $\beta \neq 1$, $\alpha \neq -2$ and so we can consider this to be equivalent to
    \[
        x(x-1) = - \frac{\alpha+1}{\alpha+2}(x-1).
    \]
    Substituting this value back into $r(x,y)$ gives
    \[
        r(x,y) = (2\alpha+1)(\alpha+1)(x-1) = 0.
    \]
    Again $\alpha \neq -\frac{1}{2}$, so we have only two cases to consider: $x = y = 1$ or $\alpha = -1$ (in which case $\beta = 0$).

    In the case $\{\alpha, \beta\} = \{-1,0\}$, we proceed to find that $V$ is a $1$-dimensional algebra. In the case $x = y = 1$, we calculate that
    \[
        (a_0a_1)(a_0a_1) - [(a_0a_1)(a_0a_1)]^\sigma = 6(2\alpha + 3)(2\alpha - 1)(a_0 - a_1).
    \]
    Thus, either we find ourselves once again in the $1$-dimensional case, or $\alpha \in \{-\frac{3}{2}, \frac{1}{2}\}$. Substituting these values of $\alpha$, $\beta$ and $x$ into (\ref{eq:C3graded}) gives the algebra products as required.
\end{proof}

\begin{prop}
    \label{prop:C3gradedideals}
    Let $V$ be one of the two $3$-dimensional algebras given in Theorem \ref{thm:C3graded}. The algebra $V$ has two non-trivial proper ideals; $\langle v_\alpha \rangle$ (or $\langle v_\beta \rangle$) if $\alpha = \frac{1}{2}$ (or $\beta = \frac{1}{2}$ respectively) and $\langle v_\alpha, v_\beta \rangle$.
\end{prop}

\begin{proof}
    Assume that we are in the case $\beta = \frac{1}{2}$ and
    \[
        (a_0a_1)(a_0a_1) = (\frac{1}{4} - \alpha)(a_0 + a_1) + (2\alpha +  \frac{1}{2})a_0a_1.
    \]
    The other case (with $\alpha = \frac{1}{2}$) is identical.

    We follow a similar method to that of Proposition \ref{prop:threeevalsideals}. We first suppose that $a_0 \in I$ and $\alpha = 0$ then, as $\beta = \frac{1}{2}$, we can conclude that $v_\beta \in I$. Then also $v_\beta v_\beta \in I$. By taking a linear combination of $a_0$, $v_\beta$ and $v_\beta v_\beta$, we conclude that $a_1 \in I$. If $a_0, a_1 \in I$ then we must have $I = V$, a contradiction and so this case does not occur.

    We now note that for any value of $\alpha$, $v_\alpha v_\alpha = v_\alpha v_\beta = 0$. In particular, this means that $\langle v_\alpha \rangle$ and $\langle v_\alpha, v_\beta \rangle$ are ideals of $V$ for all values of $\alpha$. If $v_\beta \in I$ then also
    \[
    v_\beta v_\beta = -(\alpha - \frac{1}{2})^2 v_\alpha \in I.
    \]
    As $\beta = \frac{1}{2}$. $\alpha \neq \frac{1}{2}$ and so also $v_\alpha \in I$.

    Finally, we consider the algebras $V/I$. If $I = \langle v_\alpha \rangle$ then $V/I$ is the $2$-dimensional algebra $2B(\frac{1}{2})$. If $I = \langle v_\alpha, v_\beta \rangle$ then $V/I$ is the $1$-dimensional algebra $1A$.
\end{proof}

\section{Main results}
\label{sec:main_results}

\subsection{Graded primitive axial algebras}

\begin{thm}
    Let $V = \langle \langle a_0, a_1 \rangle \rangle$ be a primitive $(\mathcal{F}, *)$-axial algebra over $\mathbb{C}$ where $|\mathcal{F}| = 2$.
    Then $V$ is either the $1$-dimensional algebra $1A$ or the $2$-dimensional algebra $2B(\alpha)$ for some $\alpha \in \mathbb{C}$. If additionally $V$ is assumed to be graded then either $V$ is the algebra $1A$ or $\alpha \in \{-1, \frac{1}{2}\}$.
\end{thm}

\begin{proof}
    As $V$ is primitive, $1$ must be a unit in the fusion law and so $(\mathcal{F}, *)$ must be a sublaw of the fusion law given in Theorem \ref{thm:2evals}. The algebras that satisfy this fusion law and their quotients are classified in Theorem \ref{thm:2evals} and Proposition \ref{prop:twoevalsideals}. The graded case follows from Corollary \ref{coro:2evalsgraded} and Proposition \ref{prop:twoevalsideals}.
\end{proof}

\begin{thm}
    Let $V = \langle \langle a_0, a_1 \rangle \rangle$ be a graded primitive $(\mathcal{F}, *)$-axial algebra over $\mathbb{C}$ where $|\mathcal{F}| = 3$. Then either $V$ is the $1A$ algebra or $V$ is isomorphic to one of the algebras given in Table \ref{tab:mainthm} where $\alpha, \beta, x \in \mathbb{C}$.
\end{thm}

\begin{landscape}
    \begin{table}
    \begin{center}
    \vspace{0.35cm}
    \noindent
    \begin{tabular}{>{$}c<{$}|>{$}l<{$}|>{$}c<{$}}
     \textrm{Basis} & \textrm{Products } & \textrm{Fusion law} \\ \hline
&& \\
     a_0,a_1 & a_0a_1 = \alpha(a_0 + a_1) &
     \begin{tabular}{>{$}c<{$}|>{$}c<{$}>{$}c<{$}}
       & 1 & \alpha  \\ \hline
     1 & 1 & \alpha \\
     \alpha & \alpha  & 1, \alpha \\
     \end{tabular} \\
     && \\
     a_0,a_1, a_0a_1 & \ml{l}{ a_0(a_0a_1) = \frac{1}{2}(1 - \beta)(\alpha + \beta)a_0 - \alpha\beta a_1 + (\alpha + \beta)a_0a_1 \\ (a_0a_1)(a_0a_1) = \frac{1}{2}(1 - 3\alpha - \beta)(\alpha + \beta)(\beta(a_0 + a_1) - a_0a_1) } &
     \begin{tabular}{>{$}c<{$}|>{$}c<{$}>{$}c<{$}>{$}c<{$}}
           & 1 & \alpha & \beta  \\ \hline
         1 & 1 & \alpha & \beta \\
         \alpha & \alpha  & 1, \alpha & \beta \\
         \beta & \beta  & \beta & 1, \alpha \\
     \end{tabular}
      \\
       && \\
     a_0,a_1, a_0a_1 & \ml{l}{ a_0(a_0a_1) = \frac{1}{2}((1 - \alpha)x a_0 - \alpha a_1 + (2\alpha + 1)a_0a_1 ) \\ (a_0a_1)(a_0a_1) = \frac{1}{4}((2\alpha - 1)(\alpha -1)x - (2\alpha + 1)\alpha)(a_0 + a_1) - \frac{1}{2}((2\alpha+1)(\alpha - 1)x + (2\alpha +3)\alpha)a_0a_1} &
     \begin{tabular}{>{$}c<{$}|>{$}c<{$}>{$}c<{$}>{$}c<{$}}
           & 1 & \alpha & \frac{1}{2}  \\ \hline
         1 & 1 & \alpha & \frac{1}{2} \\
         \alpha & \alpha  & 1, \alpha & \frac{1}{2} \\
         \frac{1}{2} & \frac{1}{2}  & \frac{1}{2} & 1, \alpha \\
     \end{tabular}
      \\
       && \\
     a_0,a_1, a_0a_1 & \ml{l}{ a_0(a_0a_1) = \frac{1}{2}((1 - \alpha)a_0 - \alpha a_1 +(2\alpha + 1)a_0a_1 \\ (a_0a_1)(a_0a_1) = \frac{1}{4}(4 - \alpha)(a_0 + a_1 - 2a_0a_1)} &
     \begin{tabular}{>{$}c<{$}|>{$}c<{$}>{$}c<{$}>{$}c<{$}}
           & 1 & \alpha & \frac{1}{2}  \\ \hline
         1 & 1 & \alpha & \frac{1}{2} \\
         \alpha & \alpha  & \frac{1}{2} & 1 \\
         \frac{1}{2} & \frac{1}{2} & 1 & \alpha \\
     \end{tabular}
     \\

    \end{tabular}
    \caption{Graded $2$-generated algebras}
    \label{tab:mainthm}
    \end{center}
    \end{table}

\end{landscape}

\begin{proof}
    As $|\mathcal{F}| = 3$, $(\mathcal{F}, *)$ must be a sublaw of one of the fusion laws given in Proposition \ref{prop:threeevalsfusionlaws} with $e = 1$ and $\alpha, \beta \in \mathbb{C}$. As the fusion laws (a) and (b) of Proposition \ref{prop:threeevalsfusionlaws} are isomorphic, by Proposition \ref{prop:isomorphicfusionlaws}, we need only consider one of these laws.

    The primitive $2$-generated algebras that obey the fusion law (a) and their quotients are classified in Theorem \ref{thm:3evals} and Proposition \ref{prop:threeevalsideals}. In Theorem \ref{thm:3evalsother}, we show that a primitive $2$-generated algebra that obeys fusion law (c) must again be an quotient of the $2$-dimensional algebra given in Corollary \ref{coro:2evalsgraded}.

    Finally, the primitive $2$-generated algebras that obey the fusion law (d) and their quotients are classified in Theorem \ref{thm:C3graded} and Proposition \ref{prop:C3gradedideals}.
\end{proof}

\subsection{Generic Jordan fusion laws}

The following results are not required in the proof of our main theorems above. However, it is of interest due to its similarity to the well-studied case of axial algebras of Jordan type.

\begin{defn}
    An axial algebra of \emph{Jordan type} is an axial algebra that obeys the following fusion law.
    \begin{center}
        \def\arraystretch{1.5}
        \begin{tabular}{>{$}c<{$}|>{$}c<{$}>{$}c<{$}>{$}c<{$}}
              & 1 & 0 & \eta  \\ \hline
            1 & 1 & 0 & \eta \\
            0 & 0  & 0 & \eta \\
            \eta & \eta  & \eta & 1, 0 \\
        \end{tabular}
    \end{center}
\end{defn}

We can see that a clear generalisation of this law is to set $0$ and $\eta$ to be generic values $\alpha, \beta \in \mathbb{C}$.

\begin{thm}
    Suppose that $V= \langle \langle a_0, a_1 \rangle \rangle$ is a primitive axial algebra with the following fusion law for some $\alpha, \beta \in \mathbb{C}$.
    \begin{center}
        \def\arraystretch{1.5}
        \begin{tabular}{>{$}c<{$}|>{$}c<{$}>{$}c<{$}>{$}c<{$}}
              & 1 & \alpha & \beta  \\ \hline
            1 & 1 & \alpha & \beta \\
            \alpha & \alpha  & \alpha & \beta \\
            \beta & \beta  & \beta & 1, \alpha \\
        \end{tabular}
    \end{center}

    Then $V$ is equal to of one of the following.
    \begin{enumerate}[i)]
        \item The $3$-dimensional algebra given in Theorem \ref{thm:3evals} where $\alpha$, $\beta$ and $x$ may take the following values
        \begin{enumerate}[a)]
            \item $\alpha = 0$, $\beta = \frac{1}{2}$, $x \in \mathbb{C}$;
            \item $\alpha = 0$, $x = \frac{\beta}{2}$;
            \item $\beta = \frac{1}{2}$, $ x = 1$;
            \item $\beta = -1$, $ x = -\frac{1}{2}$;
            \item $\beta = 2 - 3\alpha$, $ x = 1$.
        \end{enumerate}
        \item The $2$-dimensional algebra $2B(\alpha)$ where $\alpha \in \{0, -\frac{1}{2}\}$.
        \item The $2$-dimensional algebra $2B(\beta)$ where $\beta \in \{-1, \frac{1}{2}\}$.
        \item The $1$-dimensional algebra $1A$.
    \end{enumerate}
\end{thm}

\begin{proof}
    From Theorem \ref{thm:sublaw}, the algebra $V$ must be a quotient of one of the three algebras given in Theorem \ref{thm:3evals}. If it is a quotient of the $1$-dimensional algebra $1A$ then we are in case (iv).

    Next suppose that $V$ is a quotient of the $2$-dimensional algebra $2B(\alpha)$. Then $V_\beta^{a_0}$ and $V_\beta^{a_1}$ are $0$-dimensional and so $V$ must satisfy the fusion law in Corollary \ref{coro:jordansubtable}. From Corollary \ref{coro:jordansubtable}, we see that either $\alpha \in \{0, \frac{1}{2}\}$ and we are in case (ii), or $V$ is $1$-dimensional and we are in case (iv).

    Now suppose that $V$ is a quotient of the $3$-dimensional algebra given in Theorem \ref{thm:3evals}. In this case, we must have either $\beta = \frac{1}{2}$ or $x = \frac{\alpha + \beta }{2(1-\alpha)}$. The eigenspace $V_\alpha^{a_0}$ is at most one dimensional, spanned by $v_\alpha$.

    The fusion law of this algebra implies that $a_0(v_\alpha v_\alpha) = \alpha v_\alpha v_\alpha$. We can calculate that
    \[
        a_0(v_\alpha v_\alpha) - \alpha v_\alpha v_\alpha = \alpha(\alpha - \beta)(\alpha - 1)((1-\beta)x - \beta)(x-1) a_0 = 0.
    \]

    If $\alpha = 0$ then this expression holds and so the algebra is as given in Theorem \ref{thm:3evals}. If additionally $\beta = \frac{1}{2}$ then $x$ may take any value in $\mathbb{C}$. Otherwise, $x = \frac{\beta}{2}$.

    Now we can assume that $\alpha \neq 0$. Then, as $\alpha \neq 1$, $\beta \neq 1$ and $\alpha \neq \beta$, we can conclude that $((1-\beta)x - \beta)(x-1) = 0$.

    If $\beta = \frac{1}{2}$ then this becomes $\frac{1}{2}(x-1)^2$ and so we conclude that $x = 1$. Otherwise, if $\beta \neq \frac{1}{2}$ then $x = \frac{\alpha + \beta }{2(1-\alpha)}$ and we have that $(\alpha - \beta)(\beta +1)(\beta + 3\alpha -2) = 0$. Thus in this case, either $\beta = -1$ and $x = -\frac{1}{2}$, or $\beta = 2 - 3\alpha$ and $x = 1$.

    By comparing the values here with the ideals found in Proposition \ref{prop:threeevalsideals}, we see that the $2$-dimensional algebra $2B(\beta)$ can also occur, but only for the values $\beta = -1$ and $\beta = \frac{1}{2}$.
\end{proof}

\section*{Acknowledgements} I would like to thank Prof. Sergey Shpectorov for his guidance and support throughout this project.

\end{document}